\documentclass[11pt]{article}

\usepackage{amssymb,amsmath}
\usepackage{amsthm}

\newcommand{\PSbox}[3]{\framebox{\rule{0in}{#3}\includegraphics{#1}\hspace{#2}}}

\numberwithin{equation}{section}

\theoremstyle{plain}
\newtheorem{thm}{Theorem}[section]
\newtheorem{lem}[thm]{Lemma}
\newtheorem{prop}[thm]{Proposition}
\newtheorem{cor}[thm]{Corollary}

\theoremstyle{definition}
\newtheorem{defn}{Definition}[section]

\newtheorem{exmp}{Example}[section]
\newtheorem{cond}{Condition}[section]

\theoremstyle{remark}
\newtheorem*{rem}{Remark}

\author{Sam Kaufman}
\title{Delzant-type classification of near-symplectic toric 4-manifolds}
\begin{document}
\maketitle
\begin{abstract}
 Delzant's theorem for symplectic toric manifolds says that there is a one-to-one correspondence between certain convex polytopes in $\mathbb{R}^n$ and symplectic toric $2n$-manifolds, realized by the image of the moment map.  I review proofs of this theorem and the convexity theorem of Atiyah-Guillemin-Sternberg on which it relies.  Then, I describe Honda's results on the local structure of near-symplectic 4-manifolds, and inspired by recent work of Gay-Symington, I describe a generalization of Delzant's theorem to near-symplectic toric 4-manifolds. One interesting feature of the generalization is the failure of convexity, which I discuss in detail.  The first three chapters are primarily expository, duplicate material found elsewhere, and may be skipped by anyone familiar with the material, but are included for completeness.

\end{abstract}
\tableofcontents
\listoffigures

\pagestyle{plain}

\section{Background}
\subsection{Hamiltonian group actions}
Let $(M, \omega)$ be a symplectic manifold, let $G$ be a connected compact Lie group, and let $\varphi : G \times M \to M$ be a smooth Lie group action $(g,m) \mapsto g \cdot m = \varphi_{g}(m)$. $\varphi$ is \emph{symplectic} if $\varphi_g^* \omega = \omega,  \forall g \in G$.  For $\xi \in \mathfrak{g}$, we define the vector field $X_\xi$ by $X_\xi(m) = \frac{d}{dt} \vert_{t=0} (exp(t \xi) \cdot m)$.  

Let $H \in C^\infty(M)$. Define the vector field $X_H$ by $\iota(X_H)\omega = dH$. $X_H$ is called the \emph{Hamiltonian vector field associated to $H$}.  $\omega$ defines a Poisson structure on $M$ by $\{H,F\} = \omega(X_H,X_F)$. We say the group action $\varphi$ is \emph{Hamiltonian} if there exists a map $\Phi: M \to \mathfrak{g}^*$ such that for all $\xi \in \mathfrak{g}$, $X_\xi = X_{\Phi^{\xi}}$, where $\Phi^{\xi} \in C^{\infty}(M)$ is defined by $\Phi^{\xi}(m) = \langle \xi, \Phi(m) \rangle$, and furthermore, the map $j: \mathfrak{g} \to C^{\infty}(M), \xi \mapsto \Phi^{\xi}$ is a Lie algebra homomorphism.  In this case the action $\varphi$ is automatically symplectic, because $\mathcal{L}_{X_\xi}\omega = d \iota(X_\xi) \omega = d\cdot d \Phi^{\xi} = 0 $.  We call $\Phi$ the \emph{moment map} for $\varphi$.

The following two properties of the moment map will be fundamental.

\begin{prop} [Equivariance of the moment map]
\begin{equation} \label{equivariance}
Ad_g^* \Phi( g \cdot m ) = \Phi(m)
\end{equation}
\begin{proof}
\begin{equation*}
\begin{split}
\langle \xi , Ad_g^* d \Phi(g \cdot m) \rangle& = \langle g \xi g^{-1}, d \Phi (g \cdot m) \rangle \\
 & = {d \varphi_g}^* \circ \iota(X_{g \xi g^{-1}} (g \cdot m))  \omega _{g \cdot m} 
\end{split}
\end{equation*}
Now, 
\begin{equation*}
\begin{split}
X_{g \xi g^{-1}} (g \cdot m)& = \frac{d}{dt}\vert_{t=0} \exp( g (t \xi) g^{-1} ) \cdot ( g \cdot m) \\
& = \frac{d}{dt}\vert_{t=0} g \exp(t \xi)  \cdot m \\
& = {d \varphi_g} ( X_\xi (m))
\end{split}
\end{equation*}
So, 
\begin{equation*}
\begin{split}
\langle \xi , Ad_g^* d \Phi(g \cdot m) \rangle& =  \iota(X_{ \xi } ( m))  (\varphi_g^* \omega) _{ m} \\
& =   \iota(X_{ \xi } ( m))   \omega _{ m} \\
& = \langle \xi, d\Phi (m) \rangle
\end{split}
\end{equation*}

\end{proof}
\end{prop}

For the next property, we rewrite $X_{\xi} = X_{\Phi^{\xi}}$ as the requirement that for all $  v \in T_{m}M, \xi \in \mathfrak{g}$, 
\begin{equation}\label{momentadjoint}
\langle d\Phi_{m}(v), \xi \rangle = \omega( X_{\xi}(m), v)
\end{equation}
For $m \in M$, let $G_m$ be the stabilizer group of $m$ under the action $\varphi$, and let $\mathfrak{g}_m$ be its Lie algebra.
\begin{prop}[\cite{GS1}]\label{annihil}
The image of $d\Phi_m$ is $\mathfrak{g}^0_m$,the annihilator in $\mathfrak{g}^*$ of $\mathfrak{g}_m$.
\begin{proof}
The symplectic form $\omega$ defines an isomorphism $T_{m}M \to  T_{m}M^*$. Under this isomorphism, \eqref{momentadjoint} shows that the map $d\Phi:  T_{m}M \to \mathfrak{g}^*$ has the map $ \xi \mapsto \iota(X_{\xi})\omega_m, \mathfrak{g} \to T_{m}M^*$ as its transpose. The proposition then follows from linear algebra.

\end{proof}
\end{prop}

\subsection{Equivariant Moser-Darboux theorems}
In this section, we review detailed proofs of equivariant local and semi-local versions of the Darboux-Moser theorem, which we require to develop local canonical forms.  The presentation is standard and follows \cite{GS2}.

\begin{thm}[Semi-local equivariant Darboux-Moser] \label{sledm}
Let $M$ be a manifold and let $G$ be a compact connected Lie group acting on $M$. Let $X \subset M$ be a submanifold.  Let $\omega_0$ and $\omega_1$ be two $G$-invariant symplectic forms on $M$ such that $\omega_0 = \omega_1$ at $X$. Then there exists a neighborhood $U$ of $X$ and a $G$-equivariant diffeomorphism $f: U \to M$, fixing each point in $X$, such that $f^*\omega_1 = \omega_0$.
\begin{proof}
Choose a $G$-invariant Riemannian metric $\rho$ on M (e.g. by averaging over $G$, using compactness).  Let $U_0 \subset NX$ be a tubular neighborhood of the zero section of the normal bundle to $X$ such that $exp_{\rho}: U_0 \to U$ is a diffeomorphism, for $U$ a tubular neighborhood of $X$ in $M$.  Define the map $\phi_t: U \to U$ by  $\phi_t(u) = exp_{\rho}( t \cdot exp_{\rho}^{-1}(u))$, so $\phi_1 = id$, $\phi_0: U \to X$ and $\phi$ is a deformation retraction.  Define the vector field  $\xi_t = \frac{d}{dt} \phi_t$.   

Set $\omega_t = t \omega_1 + (1-t) \omega_0$.  We would like to construct a flow $f_t$ such that $f_t^* \omega_t = \omega_0$.  Set $\eta_t = \frac{d}{dt} f_t$. Setting $\sigma = \omega_1 - \omega_0$, we have
\begin{equation*}
\frac{d}{dt} ( f_t^* \omega_t ) = f_t^*( d \iota (\eta_t) \omega_t) + f_t^*(\sigma)
\end{equation*}
Now,
\begin{equation*}
\sigma - \phi_0^* \sigma = \int_0^1 \frac{d}{dt} \phi_t^* \sigma = \int_0^1 \phi_t^*(d \iota (\xi_t) \sigma) dt
\end{equation*}
But $\phi_0^* \sigma = 0$ because $\omega_1 = \omega_0$ at $X$. So
\begin{equation*}
\sigma = d \int_0^1 \phi_t^*( \iota (\xi_t) \sigma) dt
\end{equation*}
Choose $\eta_t$ such that $\iota(\eta_t) \omega_t = - \int_0^1 \phi_t^*( \iota (\xi_t) \sigma) dt$, shrinking $U_0$ as necessary so that $\omega_t$ is nondegenerate on it. Note $\eta_t = 0$ along X. Integrating $\eta_t$, we obtain a flow $f_t$ satisfying $f_t^* \omega_t = \omega_0$.  Finally, $f_t$ is G-invariant because $\eta_t$, $\xi_t$,  and $\omega_i$ are, and $\phi_t$ is G-equivariant.

\end{proof} 
\end{thm}
\begin{rem}
The global Darboux-Moser theorem presupposes a smooth family of cohomologous symplectic forms $\omega_t, t \in [0,1]$.  Here, the agreement on the submanifold $X$ and the existence of a retraction of a tubular neighbourhood onto X means the cohomology condition is trivially satisfied for the family $\omega_t = t \omega_1 + (1-t) \omega_0$.  The proof here uses a form of the Poincar\'e Lemma to construct the coboundary explicitly, and makes the $G$-equivariance explicit.
\end{rem}

One consequence of Theorem \ref{sledm} is that it allows us to linearize the $G$-action and the symplectic form simultaneously at a fixed point, as follows.

\begin{cor}\label{ptlinear} Let $(M, \omega)$ be a symplectic manifold on which $G$ acts symplectically. Let $x \in M$ be a fixed point of the $G$-action $\varphi$, so that $G$ acts on $T_xM$ by $g \mapsto d\varphi_g$. Then there exist neighborhoods $0 \in U_o \subset T_xM$ and $x \in U \subset M$ and a $G$-equivariant diffeomorphism  $h: U_o \to U$ such that $h^* \omega = \omega_x$.
\begin{proof}
As before, let $\rho$ be a $G$-invariant Riemannian metric on $M$, and let $0 \in U^{\prime}_0 \subset T_xM$ and $x \in U^{\prime} \subset M$ be such that $exp_\rho: U^{\prime}_0 \to U^{\prime}$ is a diffeomorphism. For any $v \in T_xM$, notice that
\begin{equation*}
\exp_{\rho}(d\varphi_g \cdot sv) = g \cdot \exp_{\rho}(sv)
\end{equation*}
because, by $G$-invariance of $\rho$, both are geodesics tangent to $v$ at $x$, so $\exp_{\rho}$ is $G$-equivariant.  
Set $\omega_0 = \omega_x$ and $\omega_1 = \exp_{\rho}^* \omega$ on $U_0^{\prime}$.  (Here $G$ acts by $d\varphi$.) By Theorem  \ref{sledm}, there's a $G$-equivariant diffeomorphism $f: U_0 \to U_0, f^*\omega_1 = \omega_0$, i.e. $f^* \exp_{\rho}^* \omega = \omega_0 = \omega_x$ on $U_0$, and setting $h = \exp_{\rho} \circ f : U_0 \to U$ have $ h ( g \cdot u ) = g \cdot h(u)$ as desired.

\end{proof}
\end{cor}

We can similarly linearize the $G$-action and the symplectic form along G-invariant submanifolds:

\begin{cor}\label{mfldlinear}
Let $X \subset M$ be a $G$-invariant submanifold. Then after choosing a $G$-invariant metric $\rho$, $g \mapsto d\varphi_g$ defines a $G$-action that's a bundle-map on $NX$, the normal bundle.  Then there exist neighborhoods $U_0 \subset NX$  of the zero section and $X \subset U \subset M$ and a $G$-equivariant diffeomorphism $h: U_0 \to U$ such that $h^* \omega = \omega_X$, where $\omega_X$ is any symplectic form on $NX$ that agrees with $\omega$  at the zero section (with a certain embedding $NX \to T_XM$).

\begin{proof}
Again, let $\rho$ be a $G$-invariant Riemannian metric on M, and use it to identify $NX$ with a subbundle of $T_{X}M$. Choose tubular neighborhoods $U^{\prime}_0 \subset NX$ of the zero section and $U^{\prime} \subset M$ of X such that $\exp_\rho: U^{\prime}_0 \to U^{\prime}$ is a diffeomorphism.  By $G$-invariance, $d\varphi$ maps $NX \to NX$ and $TX \to TX$.  As in Corollary \ref{ptlinear}, for $v \in N_xX$, we have
\begin{equation*}
\exp_{\rho}(d\varphi_g \cdot sv) = g \cdot \exp_{\rho}(sv)
\end{equation*}
so $\exp_\rho$ is $G$-equivariant.  Set $\omega_0 = \omega_X$ where the right-hand side is any symplectic form  which agrees with $\exp_\rho^* \omega$ along the zero section. Setting $\omega_1 = \exp_\rho^* \omega$, we can apply Theorem \ref{sledm} to conclude.  

\end{proof}
\end{cor}

\subsection{Local forms}

The following standard fact shows that the linearized actions act as complex representations.  Recall, given $(M, \omega)$, a compatible almost-complex structure is a fibre-preserving automorphism $J: TM \to TM$ such that $J^2 = -1$ and $(v,w) \mapsto \omega(v, Jw)$ defines a Riemannian metric on $M$. 
\begin{prop}
Let $(M,\omega)$ be a symplectic manifold with a compact group $G$ acting symplectically on it. Then $M$ admits a $G$-invariant compatible almost-complex structure.
\begin{proof}
Let $\rho$ be a $G$-invariant Riemannian metric on $M$.  Defining $A$ by $\rho(u,v) = \omega(u,Av)$, $-A^2$ is symmetric and positive definite.  Setting $P= \sqrt(-A^2)$ and $J=AP^{-1}$, $J$ is an invariant compatible complex structure.

\end{proof}
\end{prop}

We can use Corollary \ref{mfldlinear} to describe the neighborhood of an orbit explicitly, as follows:

\begin{lem}[Equivariant slice theorem]\label{eqslice}
Let $G$ be a compact Lie group acting (symplectically) on $M$, let $x \in M$, and let $G \cdot x$ be its orbit.  A $G$-invariant neighborhood of the orbit is equivariantly diffeomorphic to a neighborhood of the zero section in the bundle $G \times _{G_x} W$, where $G_x$ is the isotropy group of $x$ and $W = T_xM / T_x(G \cdot x)$. Furthermore, the diffeomorphism can be chosen to be a symplectomorphism when $W$ is identified with a particular embedding of $N(G \cdot x)$ in $TM$ and the symplectic form on the bundle agrees with the pullback of $\omega$ along the zero section.
\begin{proof}
As usual, choose a $G$-invariant metric on $M$ and use it to identify $N(G \cdot x)$ with a subbundle of $TM$.  $d\varphi$ is a $G$ action which is a bundle-map on $N(G \cdot x)$, and for $v \in N(G \cdot x)$, $\exp(s d\varphi_g \cdot v) = g \cdot \exp (sv)$.  Therefore, for a small neighborhood  $0 \in U_0 \subset N_x(G \cdot x)$, $\exp(U_0)$ is $G_x$-invariant, and the map
\begin{equation}
G \times N_x(G \cdot x) \to M, (g,v) \mapsto g \cdot \exp(v)
\end{equation}
descends to a map
\begin{equation}
G \times_{G_x} N_x(G \cdot x) \to M
\end{equation}
which is easily seen to be a diffeomorphism around the zero section.  The left hand side is also easily seen to be identical with the bundle $N(G\cdot x)$ and Corollary \ref{mfldlinear} gives the conclusion.

\end{proof}
\end{lem}

For an isotropic orbit of a Hamiltonian action, we can further describe the slice $W = T_xM / T_x(G \cdot x)$.

\begin{lem}\label{canbund}

Let $G$ be a compact Lie group which acts in a Hamiltonian manner on $(M,\omega)$, let $x \in M$, and let $G \cdot x$ be an isotropic orbit.  Set $V = T_x(G \cdot x)^{\omega} / (T_x(G \cdot x) \cap T_x(G \cdot x)^{\omega})$. Then $G_x$ acts on V and we can identify $W = \mathfrak{g}^0_x \times V$, where $\mathfrak{g}^0_x$ is the annihilator of $\mathfrak{g}_x$ in $\mathfrak{g}^*$. Furthermore, a neighborhood of $G \cdot x \subset M$ is equivariantly symplectomorphic to a neighborhood of zero in the bundle $Y = G \times _{G_x}( \mathfrak{g}^0_x \times V)$. Here $Y$ has the symplectic form which is the product of the one induced by the canonical form on $T^*G$ for the first two factors and some other symplectic form on $V$.
\begin{proof}
First, since $G_x$ acts symplectically and preserves $T_x(G \cdot x)$, it also preserves $T_x(G \cdot x)^{\omega}$, so $G_x$ acts on $V$ by the linear isotropy action.
This time, choose the invariant metric $\rho$ to be induced by a compatible invariant almost-complex structure $J$, i.e. 
\begin{equation}\label{compat}
\rho(v,w) = \omega(v, Jw)
\end{equation}
Then $V$ can be identified with a subspace of $N(G \cdot x)$ as follows: 
\begin{equation}
\begin{split}
T_x(G \cdot x)^{\omega} & = \{ v \in T_xM | \omega(v,w) = 0, \forall w \in T_x(G \cdot x)\} \\
& = \{ v \in T_xM | \rho(v, -Jw) = 0 , \forall w \in T_x(G \cdot x) \}
\end{split}
\end{equation}
so $T_x(G \cdot x)^{\omega} = (J T_x(G \cdot x))^{\bot} = J( T_x(G \cdot x))^{\bot}$. Since $G \cdot x$ is isotropic, \ref{compat} shows that $J(T_x(G \cdot x)) \perp T_x(G \cdot x)$.  So we can identify
\begin{equation}
V  \cong J( T_x(G \cdot x))^{\bot} \cap  T_x(G \cdot x)^{\bot} 
\end{equation}
and this gives the orthogonal splitting 
\begin{equation}
T_xM  \cong   T_x(G \cdot x) \oplus J( T_x(G \cdot x))  \oplus V 
\end{equation}
which, grouping the first two terms together, is also a symplectic splitting, because
\begin{equation}
 ( T_x(G \cdot x) \oplus J( T_x(G \cdot x)) )^{\omega} \cong V 
\end{equation}

Furthermore, for the moment map $\Phi$, from equation \eqref{momentadjoint} we have that 
\begin{equation}
\ker d\Phi_x = T_x(G\cdot x)^{\omega}
\end{equation}
By definition, $V \subset T_x(G \cdot x)^\omega$, and $T_x(G \cdot x) \subset  T_x(G \cdot x)^\omega$ because $G \cdot x$ is isotropic. So we can identify 
\begin{equation}
T_xM / \ker d\Phi_x  \cong J ( T_x(G \cdot x))
\end{equation} 
So $d\Phi_x: J(T_x(G \cdot x)) \to \mathfrak{g}^*$ is an isomorphism onto its image, and by Proposition \ref{annihil}, $d\Phi_x(T_xM) = \mathfrak{g}^0_x.$
Summarizing, we have: 
\begin{equation}
T_xM = (T_x(G \cdot x) \oplus \mathfrak{g}^0_x)\oplus V 
\end{equation}
where all direct sums are orthogonal splittings, and the second one is a symplectic splitting.

Finally, applying the equivariant slice theorem to this splitting we obtain an equivariant symplectomorphism to the bundle
\begin{equation}
Y \cong G \times_{G_x} ( \mathfrak{g}^0_x \times V)
\end{equation}
where the $G_x$ action on $( \mathfrak{g}^0_x \times V)$ is the product of the co-adjoint action on $\mathfrak{g}^0_x$ (by equivariance of the moment map) and and the linear isotropy action on V.  The natural symplectic form on this bundle is as described in the statement (we can take the constant form $\omega_x$ on the $V$ factor), and by the identification via $\Phi$ it's of the form required by the equivariant slice theorem. 

\end{proof}
\end{lem}

The remainder of this section will be devoted to analyzing the moment map in the above case.  The calculation follows \cite{LT}, Lemma 3.5. 

To see the moment map on $Y$, it's convenient to construct $Y$ via symplectic reduction of a space on which the moment map is easy to calculate.  This is done as follows.  

Recall, $T^*G$, like any cotangent bundle, has a canonical 1-form $\theta$ defined, for $p \in T^*G, v \in T_P(T^*G)$ by
\begin{equation*}
\langle \theta_p, v \rangle = \langle p, d \pi_pv \rangle
\end{equation*}
where $\pi : T^*G \to G$ is the projection, which induces a canonical symplectic form $\omega_G = -d\theta_G$.  Let $H \subset G$ be a subgroup, and let $H$ act on $G$ by the right action $h \cdot g \mapsto g h^{-1}$.   

The actions of $H$ and $G$ on $G$ (like any diffeomorphism) induce symplectomorphisms of $(T^*G, \omega_G)$:

Writing $(q,p) \in T^*G, p \in T^*_qG$, we define the action of $G$ by 
\begin{equation}
g \cdot (q,p) \mapsto (g \cdot q, (dg^{-1})^*p)
\end{equation}
but 
\begin{equation}
\begin{split}
\langle \theta_p , v \rangle & = \langle p, d\pi_p v \rangle \\
& = \langle (dg^{-1})^* p, dg d\pi_p v \rangle \\
& = \langle g^* \theta_p,  v \rangle 
\end{split}
\end{equation}

Similarly, the action of $H$ is
\begin{equation}
h \cdot (q,p) \mapsto ( q  \cdot h^{-1}, (dh_R)^*p)
\end{equation}
but 
\begin{equation}
\begin{split}
\langle \theta_p , v \rangle & = \langle p, d\pi_p v \rangle \\
& = \langle (dh_R)^* p, dh^{-1}_R d\pi_p v \rangle \\
& = \langle h^* \theta_p,  v \rangle 
\end{split}
\end{equation}

Since both actions preserve $\theta_G$, both preserve $\omega_G$.

 Since $g^*\theta_G = \theta_G$ we have, for $\xi \in \mathfrak{g}$:
\begin{equation}
0 = \mathcal(L)_{X_{\xi}} \theta_G = \iota(X_{\xi}) d \theta_G + d \iota(X_{\xi}) \theta_G
\end{equation}
so 
\begin{equation}
\iota(X_{\xi}) \omega_G = d ( \iota( X_{\xi}) \theta_G)
\end{equation}
and similarly for the right $H$ action, so both actions are Hamiltonian with moment map 
\begin{equation}
\xi \mapsto \iota(X_{\xi}) \theta_G (p) = \langle p, d\pi X_{\xi} \rangle
\end{equation}

In the case of the (left) $G$-action, we have, for $p \in T^*G, p = (g,v)$, 
\begin{equation}
d \pi_p(X_{\xi}) = (dg_R) \xi
\end{equation}
so the moment map is
\begin{equation}
p = (g,v) \mapsto(  \xi \mapsto \langle v, (dg_R) \xi \rangle  ) 
\end{equation}
or
\begin{equation}
\hat{\Phi}_G^{\xi} = \langle (dg_R)^*v, \xi \rangle
\end{equation}
Similarly, for the (right) $H$-action, we have
\begin{equation}
\hat{\Phi}_H^{\xi} = \langle -j^*(dg)^*v, \xi \rangle
\end{equation}
where $j: \mathfrak{h} \to \mathfrak{g}$ is the inclusion.

Suppose $H$ also acts in a Hamiltonian way on a symplectic vector space $(V, \omega_V)$, with moment map $\Phi^V_H: V \to \mathfrak{h}^*$.  We can define two commuting Hamiltonian actions on $T^*G \times V$ by 
\begin{subequations}
\begin{align}
g^{\prime} \cdot (g, v) & \xrightarrow{\varphi_G} (g^{\prime}g, v) \\ 
h \cdot (g, v) & \xrightarrow{\varphi_H} (gh^{-1}, h \cdot v) 
\end{align}
\end{subequations}
In what follows, let $(g, \eta, v)$ be coordinates on $T^*G \times V$, with $g \in G, \eta \in T_g^*G, v \in V$. Denoting the $G$ and $H$ -action's respective moment maps by $\Phi_G: T^*G \times V \to \mathfrak{g}$ and $\Phi_H: T^*G \times V \to \mathfrak{h}$, by the previous calculations we have
\begin{subequations}
\begin{align}
\Phi_G( g, \eta, v) &= (dg_R)^* \eta \label{PhiG}\\
\Phi_H( g, \eta, v ) &= -j^*(dg)^* \eta + \Phi^V_H(v)
\end{align}
\end{subequations}

We would now like to consider the symplectic reduction of $T^*G \times V$ with respect to the $H$ action, and show that the resulting space, $\Phi_H^{-1}(0) / H$, is our model bundle $Y$.

We can write the set $\Phi_H^{-1}(0)$ as
\begin{equation}\label{Hzero}
\Phi_H^{-1}(0) = \{ ( g, \eta, v ) | \langle -j^*(dg)^* \eta + \Phi^V_H(v), \xi \rangle =0, \forall \xi \in \mathfrak{h} \}
\end{equation}

For each fixed pair $(g,v)$, we can identify the fibre $\Phi^{-1}_H(0) |_{g \times T^*_gG \times v}$ with $\mathfrak{h}^0 \subset \mathfrak{g}^*$, the annihilator of $\mathfrak{h}$.  First, by \eqref{Hzero}, the fibre's an affine subspace of $T_g^*G$.  Let $A: \mathfrak{g} \to \mathfrak{h}$ be any $H$-equivariant projection (so $A \circ j = \text{Id}$). The map $ \mathfrak{k} \mapsto (dg^*)^{-1}( \mathfrak{k} + A^*\Phi^V_H(v))$ is then an isomorphism $ \mathfrak{h}^0 \to \Phi^{-1}_H(0) |_{g \times T^*_gG \times v}$.  Thus we've identified
\begin{equation}
G \times \mathfrak{h}^0 \times V \cong \Phi_H^{-1}(0) \subset T^*G \times V
\end{equation}
via the map
\begin{equation}
(g, \mathfrak{k}, v) \mapsto ( g, (dg^*)^{-1}( \mathfrak{k} + A^*\Phi^V_H(v)), v)
\end{equation}

In these coordinates, we can rewrite \eqref{PhiG} as 
\begin{equation}
\begin{split}
\Phi_G(g, \mathfrak{k}, v)  & = (dg_R)^* (dg^*)^{-1} (\mathfrak{k} + A^* \Phi^V_H(v) )\\
& = Ad^*(g)( \mathfrak{k} + A^* \Phi^V_H(v))
\end{split}
\end{equation}

Since the $G$ action commutes with the $H$ action, it descends to the quotient space with the same moment map.  Setting $H = G_x$, so $\mathfrak{h}^0 = \mathfrak{g}_x^0$, it's clear that the quotient $\Phi_H^{-1}(0) / H = G \times_{G_x} ( \mathfrak{g}_x^0 \times V)$ as desired.  We've now shown:

\begin{lem}[\cite{LT}, Lemma 3.5 ] \label{canmoment} 
For any $G_x$-equivariant projection $A: \mathfrak{g} \mapsto \mathfrak{g}_x^0$, there's a symplectic structure on $Y$ such that Lemma \ref{canbund} is true and the moment map on $Y$ is given by $\Phi_Y ([g, \eta, v]) = Ad^*(g)(\eta + A^*\Phi_V(v))$.
\end{lem}

\section{Hamiltonian torus actions}
In this section, we apply the canonical local forms to the case where the group $G$ is a torus $T^n$.
\subsection{Convexity}
\begin{prop}
Let $G = T^n$. Then $(1)$ $G$-orbits are isotropic, and $(2)$ $\Phi$ is constant on $G$ orbits.
\begin{proof}
The same calculation shows both. For $\eta, \xi \in \mathfrak{g}$:
\begin{equation*}
\begin{split}
\iota(X _{\eta})(d\Phi^{\xi}) &= \omega( X_{\eta}, X_{\xi}) \\
&= \{\Phi^{\eta}, \Phi^{\xi}\} \\
&= j([\eta,\xi]) = 0
\end{split}
\end{equation*}
( Claim (2) also follows directly by equivariance of $\Phi$ since $Ad^{T^n} = Id$ .)

\end{proof}
\end{prop}

The following lemma is quoted without proof from elementary representation theory.
\begin{lem}\label{weights}
Let $T^n$ act linearly and unitarily on $\mathbb{C}^m$.  Then there exists an orthogonal decomposition $\mathbb{C}^m = \oplus_{k=1}^m V_{\lambda^{(k)}}$ into one-dimensional $T$-invariant complex subspaces and linear maps $\lambda^{(i)} \in \mathfrak{t}^*, i=1,...,m$ such that on $V_{\lambda^{k}}$, $T^n$ acts by $(e^{it_1},...,e^{it_n}) \cdot v = e^{i \sum_j \lambda_j^{(k)}t_j}v$.
\end{lem}
The covectors $\lambda^{(k)}$ are called the \emph{weights} of the representation.

\begin{cor}[Local convexity \cite{GS1}]\label{locconv}
Let $x \in M$ be a fixed point of a Hamiltonian $T^n$ action with moment map $\Phi$, and let $p=\Phi(x)$. Then there exist open neighborhoods $x \in U \subset M $ and $p \in U^{\prime} \subset \mathfrak{t}^*$ such that $\Phi(U) = U^{\prime} \cap ( p + S( \lambda^{(1)}, ..., \lambda^{(k)}))$, where $\lambda^{(i)}$ are the weights of the isotropy representation $g \mapsto dg$ on $T_xM$ and $S(   \lambda^{(1)}, ..., \lambda^{(k)}) = \{ \sum_{i=1}^n s_i \lambda^{(i)}, s_i \geq 0\}$.
\begin{proof}
By choosing a compatible invariant almost-complex structure $J$ on $M$ and the induced invariant Riemannian metric $\rho$, we make $T_xM$ into a complex vector space and obtain a unitary representation $g \mapsto dg$. Furthermore, any one-dimensional complex subspace is symplectic, since $\omega_x(v,Jv) = \rho(v,v) \neq 0$.  Thus the $T$-invariant subspaces $V_{\lambda}$ are symplectic, and they are pairwise symplectically orthogonal because they're $J$-invariant and $\rho$-orthogonal. We can therefore write $\omega_x = \sum_{i=1}^m dz_i \wedge d\bar{z_i} = \omega_0$, where $z_i$ is a complex coordinate on $V_{\lambda^{(i)}}$, and this identifies $(T_xM,\omega) \cong (\mathbb{C}^m, \omega_0)$.  The induced $T^n$ action on $(\mathbb{C}^m, \omega_0)$ is described in Lemma \ref{weights} and its moment map is $z \mapsto \sum_i |z_i|^2 \lambda^{(i)}$, which can be checked easily. Finally, by the equivariant Darboux theorem, a neighborhood $0 \in \hat{U} \subset (\mathbb{C}^m,\omega_0)$ is equivariantly symplectomorphic to a neighborhood $x \in U \subset M$, so the image of the moment map is the same, up to translation.

\end{proof}
\end{cor}

\begin{cor}[Relative local convexity \cite{GS1}]\label{rellocconv}
Let $x \in M$ have orbit $T \cdot x$ and let $T_x$ be the isotropy group of $x$.  Let $p = \Phi(x)$ and let $\lambda^{(1)}, ..., \lambda^{(k)}$ be the weights of the isotropy representation on a slice $V$ at $x$.  Let $\pi: \mathfrak{t}_x \to \mathfrak{t}$ be the inclusion. Then there exist neighborhoods $U \subset M$ of $T \cdot x$ and $U^{\prime}$ of p such that $\Phi(U) = U^{\prime} \cap ( p + S^{\prime}( \lambda^{(1)}, ..., \lambda^{(k)}))$ where $S^{\prime}( \lambda^{(1)}, ..., \lambda^{(k)}) = (\pi^*)^{-1}S(\lambda^{(1)}, ..., \lambda^{(k)})$ .
\begin{proof}
By Lemma \ref{canmoment}, around $T \cdot x$, we have $\Phi ([g, \eta, v]) = Ad^*(g)(\eta + A^*\Phi_V(v))$, up to translation. Since $G=T, \text{Ad} = \text{Id}$, so $\Phi ([g, \eta, v]) = p + \eta + A^*\Phi_V(v)$. By the previous corollary, $\Phi_V(V) = S(\lambda^{(1)}, ..., \lambda^{(k)}) \subset \mathfrak{t}_x^*$. Finally, for any projection $A$, the set $\{ \eta + A^*S(\lambda^{(1)}, ..., \lambda^{(k)}) | \eta \in \mathfrak{t}_x^0\}$ is equal to $S^{\prime}( \lambda^{(1)}, ..., \lambda^{(k)})$.

\end{proof}
\end{cor}

Recall, a function $f: M \to \mathbb{R}$ is Bott-Morse if each component of its critical set $C_f$ is a submanifold of $M$, and for each $x \in C_f$, the Hessian $d^2f_x$ is nondegenerate on $N_xC_f$. (The index of $d^2f_x$ is constant on each component of $C_f$.) The following lemma is key:

\begin{lem}[\cite{GS1}]
For each $\xi \in \mathfrak{g}$, $\Phi^{\xi}$ is Bott-Morse, and the indices and coindices of its critical manifolds are all even.
\begin{proof}
From our canonical local form, 
\begin{equation*}
\Phi^{\xi}([g, \eta, v]) = \langle p, \xi \rangle + \langle \eta, \xi \rangle + \sum|z_i|^2 \langle A^*\lambda^{(i)}, \xi \rangle
\end{equation*}
Modulo some mess from the quotient, the result can be read off: for $x \in C_{\Phi^{\xi}}$ have $\xi \in \mathfrak{g}_x$, so $\eta$ is free to vary in $\mathfrak{g}_x^0$.  We also have there that $z_i = 0$ or $z_i$ is free and $\langle A^*\lambda^{(i)}, \xi \rangle = 0$, so the critical sets are manifolds. The index is $2k$ where $k$ is the number of $i$'s such that $\langle A^*\lambda^{(i)}, \xi \rangle < 0$, and similarly for coindex, because each $V_{\lambda}$ is 2-dimensional.

\end{proof}
\end{lem}

\begin{lem}[\cite{GS1}]\label{locmax}
For each $\xi \in \mathfrak{g}$, $\Phi^{\xi}$ has a unique connected component of local maxima.
\begin{proof}
Let $C_1,...,C_k$ be the connected critical manifolds of $\Phi^{\xi}$ consisting of local maxima, and let $C_{k+1},...,C_N$ be the remaining connected critical manifolds. $M = \amalg_{i=1}^N W_i$, where $W_i$ is the stable manifold of $C_i$. Note $\text{dim}(W_i) = \text{index}(C_i) + \text{dim}(C_i)$. For $i=1,...,k$, $\text{dim}(W_i) = \text{dim}(M)$, so $W_i, i=1,...,k$ is open. For $i=k+1,...,N$, codim $C_i \geq 2$.  Therefore $M \setminus \cup_{i=k+1}^N W_i$ is connected, i.e. $\cup_{i=1}^k W_i$ is connected, so $k=1$, and there is a unique connected component of local maxima.

\end{proof}
\end{lem}

\begin{cor}[Global convexity \cite{GS1}]
$\Phi(M) \subset \mathfrak{t}^*$ is a convex polytope, specifically the convex hull of the image of the fixed points, $\Phi(M^T)$.
\begin{proof}
Let $p \in \partial \Phi(M), x \in \Phi^{-1}(p)$.  By Corollary \ref{rellocconv}, there exist neighborhoods  $U \subset M$ of $x$ and $U^{\prime}$ of p such that $\Phi(U) = U^{\prime} \cap ( p + S^{\prime}( \lambda^{(1)}, ..., \lambda^{(k)}))$ where $\lambda^{(1)}, ..., \lambda^{(k)}$ are the weights of the isotropy representation on a slice $V$ at $x$. We can choose $\xi \in \mathfrak{t}$ such that $\langle \cdot , \xi \rangle = 0$ on a boundary component of $S^{\prime}$ and $\langle \cdot , \xi \rangle < 0$ on $S^{\prime}$. Then if $\langle \Phi(p) , \xi \rangle = a$, $\langle \Phi(x) , \xi \rangle \leq a$ for $x \in U$, i.e. $a$ is a local maximum of $\Phi^{\xi}$, so by Lemma \ref{locmax}, $\Phi^{\xi} \leq a$ on $M$.  Repeating this argument for each face of $S^{\prime}$, we have $\Phi(M) \subset p + S^{\prime}(\lambda^{(1)}, ..., \lambda^{(k)} )$. Applying this argument to all boundary components of $\Phi(M)$, $\Phi(M)$ is convex.  Finally, by the local canonical form, if $\Phi(x)$ is an extremal point of $\Phi(M)$, we must have $\mathfrak{g}_x^0 = \emptyset$, so $x$ is a fixed point.

\end{proof}
\end{cor}

The following connectedness result requires a more involved Morse-theoretic argument, and will be quoted without proof.
\begin{lem}[Connectedness \cite{At}, \cite{LT}]\label{connected}
For every $a \in \mathfrak{t}^*$, the fiber $\Phi^{-1}(a)$ is connected.
\end{lem}
\subsection{Delzant's theorem}
In this section we consider the case of an effective Hamiltonian $T^n$ action on $M^{2n}$ (effective means that the action has trivial kernel.)  In this case we say $M$ is a \emph{toric} $2n$-manifold.  In what follows we'll often write $\Delta = \Phi(M)$.

\begin{prop}[Smoothness]\label{Smoothness}
Let $T^m$ act linearly on $(\mathbb{C}^n, \omega_0)$. By Lemma \ref{weights}, we have
$\mathbb{C}^n = \oplus_{k=1}^n V_{\lambda^{(k)}}$ such that on $V_{\lambda^{k}}$, $T^m$ acts by $(e^{it_1},...,e^{it_n}) \cdot v = e^{i \sum_j \lambda_j^{(k)}t_j}v$, i.e. the action factors through a map $T^m \to ^{\Psi} T^n, \exp(t) \mapsto \exp( \langle \lambda^{(1)}, t \rangle, ..., \langle \lambda^{(n)}, t \rangle)$. (So $\lambda^{(k)} \in \mathbb{Z}^m$.)Then if $T^m$ acts effectively, $m \leq n$.  If $m = n$, $\lambda^{(k)}$ are a $\mathbb{Z}$-basis of $\mathbb{Z}^m \cong \mathfrak{t}^*$.

\begin{proof}
The map $T^m \to^{\Psi} T^n$ lifts to the linear map $\mathfrak{t}^m \to^{\psi} \mathfrak{t}^n$ given by the weights.  If $m > n$, $\psi$ and hence $\Psi$ has nontrivial kernel, contradicting effectiveness.  Similarly, if $m=n$, $\psi$ must have trivial kernel, i.e. be an isomorphism.  In this case, if $\{\lambda^{(k)}\}$ is not a $\mathbb{Z}$-basis of $\mathbb{Z}^m$, then there exist lattice points in $\mathfrak{t}^n$ that are not the images of lattice points in $\mathfrak{t}^m$. Since $\psi$ is onto, this means that $\Psi$ has nontrivial kernel, contradicting effectiveness.

\end{proof}
\end{prop}

\begin{cor}\label{delzantcond}
For an effective Hamiltonian $T^n$ action in $M^{2n}$, the moment polytope $\Delta = \Phi(M)$ satisfies the following \emph{Delzant conditions}: (1) \emph{simplicity} - n edges meet at each vertex (2) \emph{rationality} - each vertex is of the form $\{p + \sum t_i v_i, t_i \geq 0, v_i \in \mathfrak{t}^*\}$, such that $v_i$ has integral entries (3) \emph{smoothness} - at each vertex, $\{v_i\}$ is a $\mathbb{Z}$-basis of $\mathbb{Z}^n$.
\begin{proof}
Apply Proposition \ref{Smoothness} and the analysis in the proof of Corollary \ref{locconv} to the fixed points of the action.

\end{proof}
\end{cor}

\begin{prop}\label{momentbijection}
The map $M/T \to \Delta, m \mapsto \Phi(m)$ is a bijection.
\begin{proof}
Let $x \in M, p=\Phi(x)$. Recall our local model for a neighborhood of $T \cdot x \subset M$, i.e. a neighborhood of the zero section in the bundle $T \times _{T_x} (\mathfrak{t}_x^0 \times V)$ with moment map $[g, \eta, v] \mapsto p + \eta + \sum |z_i|^2 A^* \lambda_{(i)}$, with $z_i$ coordinates on $V$. By definition of $A$ and $\mathfrak{t}_x^0$, $\text{image}( A^*) \cap \mathfrak{t}_x^0 = \emptyset$, so if the set $\{A^* \lambda_{(i)}\}$ is independent, then $\Phi: M/T \to \Delta$ is locally a bijection onto its image.  Since $A^*$ is injective, it's sufficient to show that the weights are independent.  Let dim $T_x = k$. Then dim $V=2k$.  The action of $T_x$ on $V$ is also effective, because otherwise, by the local form the action of $T$ wouldn't be effective. So by Proposition $\ref{Smoothness}$, $\lambda_{(i)}$ are independent.  Therefore $\Phi: M/T \to \Delta$ is locally a bijection.  To see that it's a global bijection,  use Lemma \ref{connected} to see that each set $\Phi^{-1}(a)$ must be a single orbit.

\end{proof}
\end{prop}

Any polytope $\Delta \subset \mathfrak{t}^*$ satisfying the conditions of Corollary \ref{delzantcond} is called a \emph{Delzant polytope}. Delzant \cite{De} proved that for any Delzant polytope $\Delta$, there exists a unique symplectic toric manifold $M$ such that $\Phi(M) = \Delta$.  The construction for the existence proof can be found in essentially the same form in any of \cite{De, G, LT, dS}, and will be skipped. The uniqueness proof I reproduce below is due to \cite{LT}.

\begin{thm}[Delzant \cite{De}]\label{delzant}
Let $M_1, M_2$ be two compact, connected, symplectic toric $2n$-manifolds with moment maps $\Phi_1, \Phi_2$ and moment polytopes $\Delta_1, \Delta_2$. If $\Delta_1 = \Delta_2$, then there exists a $T^n$-equivariant symplectomorphism $f: M_1 \to M_2$ such that $\Phi_2 \circ f = \Phi_1$.  
\end{thm}

The proof involves several intermediate results.

\begin{prop}\label{deltadeterm}
Let $\alpha \in \Delta$. Then a neighborhood of $\Phi^{-1}(\alpha)$ is determined by $(\Delta, \alpha)$.
\begin{proof}
By the canonical local form, we need to determine the subspace $\mathfrak{g}_x^0$ and the weights $\lambda^{(i)} \in \mathfrak{g}_x^*$.  First, note that the subspace $\mathfrak{g}_x^0 \subset \mathfrak{g}^*$ is the subspace parallel to the affine face of $\Delta$ that $\alpha$ belongs to.  Let $\mathfrak{g}_x^0$ have codimension $k$.  By simplicity, each codimension $k$ face $F_k$ belongs to $k$ codimension $k-1$ faces $\{F_{k-1}^i\}_{i=1}^k$. For each codimension $k-1$ face $F_{k-1}^i$, we can choose a covector $\alpha_{k,i} \in F_{k-1}^i \setminus F_k$, eg $\alpha_{k,i} \in  F_{k-1}^i \cap F_k^{\bot}$.  The covector $\alpha_{k,i} $ is the image $A^* \lambda^{(i)}$ for some projection $A$ corresponding to the choice of $\alpha_{k,i}$, but this choice is irrelevant because any two such models are symplectomorphic.

\end{proof}
\end{prop}

\begin{cor}\label{locsymp}
Suppose $\Delta = \Delta_1 = \Delta_2$. Then for any $\alpha \in \Delta$, there exists a neighborhood $U$ of $\alpha$ such that, for ${M_1}_U = \Phi_1^{-1}(U)$ and ${M_2}_U = \Phi_2^{-1}(U)$, there exists a $T^n$-equivariant symplectomorphism $f: {M_1}_U \to {M_2}_U$ such that $\Phi_2 \circ f = \Phi_1$.
\begin{proof}
By Proposition \ref{deltadeterm}, both ${M_1}_U$ and  ${M_2}_U$ are equivariantly symplectomorphic to the same canonical model.

\end{proof}
\end{cor}

To piece these local symplectomorphisms together, we follow $\cite{LT}$ and use sheaf cohomology.

Let $\mathcal{U}$ be a cover of $\Delta$ with the property that each $U \in \mathcal{U}$ has the property in Corollary $\ref{locsymp}$, and let $\mathcal{H}_U$ be the set of all moment-preserving $T$-equivariant symplectomorphisms ${M_1}_U \to {M_1}_U$.  Notice that this group is abelian: $\mathcal{H}_U$ acts on fibers, the fibres are $G/G_x$, and the only $G$-equivariant diffeomorphisms $G/G_x \to G/G_x$ are multiplications by elements of $G$.  Since $G$ is commutative,   $\mathcal{H}_U$ is commutative on fibres, so commutative.  We can therefore use the groups $\mathcal{H}_U$ to define a sheaf of abelian groups and its sheaf cohomology $H^*(\Delta, \mathcal{H})$.

By Corollary $\ref{locsymp}$, $\mathcal{U}$ defines a 1-cochain in the sheaf as follows. For each $U_i \in \mathcal{U}$ choose a $T$-equivariant moment preserving symplectomorphism $f_i: {M_1}_{U_i} \to {M_2}_{U_i}$. For each pair $U_i, U_j \in \mathcal{U}$, set $h_{ij} = f_i^{-1} \circ f_j \in \mathcal{H}_{U_i \cap U_j}$. This is a cocycle: for $U_i, U_j, U_k \in \mathcal{U}, U_i \cap U_j \cap U_k \neq \emptyset$, $h_{ij} \circ h_{jk} \circ h_{ki} = f_i^{-1} \circ f_j \circ f_j^{-1} \circ f_k \circ f_k^{-1} \circ f_i = Id$, so $h_{ij}$ defines a cohomology class in $H^1(\Delta, \mathcal{H})$.

Suppose $h_{ij}$ is a coboundary, i.e. $h_{ij} = h_i \circ h_j^{-1}$, for $h_i, h_j \in \mathcal{H}_{U_i}$.
Then $f_i^{-1} \circ f_j = h_i \circ h_j^{-1}$, or $f_j \circ h_j = f_i \circ h_i$ on $U_i \cap U_j$.  Then the map $x \in {M_1}_{U_i} \mapsto (f_i \circ h_i)(x)$ is a globally well-defined $T$-equivariant moment preserving symplectomorphism $M_1 \to M_2$.

So, we will show that $H^1(\Delta, \mathcal{H}) = 0$ by showing that $H^k(\Delta, \mathcal{H}) = 0, \forall k > 0$.

Define auxiliary sheafs as follows:

Let $\underline{\ell \times \mathbb{R}}$ be the locally constant sheaf on $\Delta$ with values in the abelian group $\ell \times \mathbb{R}$, where $\ell$ is the integer lattice $\mathbb{Z}^n \subset \mathfrak{t}$.

For each $U \in \Delta$, let $\tilde{C^{\infty}}(U)$ be the set of smooth $T$-invariant functions on $M_U$, so $f \in \tilde{C^{\infty}}(U) \implies f = h \circ \Phi$ for some smooth function $h$ on $U$.  Call this sheaf $C^{\infty}$.

Define a map $j: \underline{\ell \times \mathbb{R}} \to C^{\infty}$ by $j(\xi,c)(x) = \langle \xi , \Phi(x) \rangle + c, \forall x \in M_U$.

Define a map $\Lambda: C^{\infty} \to \text{Symp}(M_U)$ by $\Lambda(f)(x) = \exp(X_f)$ where $X_f$ is the Hamiltonian vector field associated to $f$ (this is the time-1 flow).

\begin{prop}\label{lambdatoh}

 $\Lambda: C^{\infty} \to \mathcal{H}$, i.e. the Hamiltonian flow preserves $\Phi$ and is $T$-equivariant.
\begin{proof}
$\exp(X_f)$ preserves $\Phi$: For $\xi \in \mathfrak{t}$,
\begin{equation*}
\iota(X_f) d\Phi^{\xi} = \omega(X_f,X_{\xi}) = -\iota(X_{\xi})df = 0
\end{equation*}
by $T$-invariance of f.

$\exp(X_f)$ is $T$-equivariant:  Have $f(x) = f(t \cdot x) \implies df_x = df_{t\cdot x} \circ dt_x $. Since $T$ acts symplectically, have
\begin{equation*}
\iota(dt_x {X_f}_x) \omega_{t\cdot x} \circ dt_x = \iota({X_f}_x) \omega_x = df_x =  df_{t\cdot x} \circ dt_x
\end{equation*}
Cancelling $dt_x$, have
\begin{equation*}
\iota(dt_x {X_f}_x) \omega_{t\cdot x} =  df_{t\cdot x}
\end{equation*}
i.e. $dt \circ X_f = X_f$, so $X_f$ is $T$-equivariant, and its flow is also.

\end{proof}
\end{prop}

\begin{lem}
The sequence of sheaves $0 \to \underline{\ell \times \mathbb{R}} \to^{j} C^{\infty} \to^{\Lambda} \mathcal{H}$ is exact.
\begin{proof}
$j$ is injective because any open set in $\Delta$ suffices to determine an affine function.

$im(j) \subset ker( \Lambda)$ since, for any $(\xi, r) \in \ell \times \mathbb{R}$, $j((\xi,r))(x) = \Phi^{\xi}(x) + r$, so it's the moment for $\xi$. By definition, $\exp X_{\Phi^{\xi}}(x) = \exp(\xi) \cdot x = id \cdot x$ because $\xi \in \ell$.

$ker( \Lambda) \subset im(j)$: Let $f \subset \tilde{C^{\infty}}(U), \Lambda(f) = id, f=h \circ \Phi$.  Since the flow of f is $G$-invariant and tangent to the orbits, at each point $x \in \Delta$, there exists $ \xi_x \in \mathfrak{g}$ such that on $\Phi^{-1}(x)$, $X_f = X_{\xi_x}$, and locally $\xi_x$ can be chosen to be continuous.  On the interior of $\Delta$, the $G$-action is free. So on $\text{int}(\Delta)$,  $Id = \exp(X_f) = \exp(X_\xi) = \exp(\xi) \cdot x \implies \exp(\xi) = id \in G$, or that $\xi_x \in \ell$ for $ x \in \text{int}(\Delta)$.  Since $\ell$ is discrete and $\xi_x$ is continuous, we must have that $\xi_x$ is locally constant on $\text{int}(\Delta)$, so also on $\partial \Delta$. So $df_x = d\Phi^{\xi}_x = d \langle \Phi(x), \xi \rangle \implies f(x) = \langle \Phi(x), \xi \rangle +r$, as claimed.

\end{proof}
\end{lem}

Will now show that $\Lambda: C^{\infty} \to \mathcal{H}$ is surjective.

Choose $\alpha \in U \subset \Delta$ simply-connected such that $M_\alpha$ is a deformation retract of $M_U$ (this is possible by our local model.) Let $f: M_U \to M_U \in \mathcal{H}(U)$.  By $G$-equivariance, can write $f(p) = \gamma( \Phi(p)) \cdot p$ where $\gamma: U \to G$ is smooth.  Since $\pi_1(U) = 0$, can lift $\gamma$ to a map $\tilde{\gamma}: U \to \mathfrak{g}$ such that $\gamma = \exp \tilde{\gamma}$.  So have $f(p) = \gamma( \Phi(p)) \cdot p = (\exp (\tilde{\gamma} (\Phi(p)))) \cdot p = \exp( X_{ \tilde{\gamma}(\Phi(p))} (p))$, i.e. $f$ is the time-1 flow of the vector field $p \mapsto X_{\tilde{\gamma}(\Phi(p))}$.  We would like to show that $Y = X_{\tilde{\gamma} \circ \Phi}$ is a Hamiltonian vector field.  

Set $U_0 = U \cap \text{int}(\Delta)$.  As we've seen, $M_{U_0}$ is a principal $G$-bundle.
\begin{prop}
$\iota(Y) \omega |_{M_{U_0}}$ is a basic form on this bundle.
\begin{proof}
We need to show that (1) $\iota(Y) \omega$ is $G$-invariant, and (2) that for all vectors $v$ tangent to the fibre, $\iota(v) \iota(Y) \omega = 0$.

(1):  Let $v \in T_xM$. Then
\begin{equation*}
\begin{split}
g^*( \iota(Y) \omega )_x (v) &= ( \iota(Y) \omega)_{g \cdot x} ( dg_x v) \\
&= \iota( {X_{ \tilde{\gamma} \circ \Phi ( g \cdot x) }}_{g \cdot x}) \omega_{g \cdot x} ( dg_x v) \\
&= \iota( {X_{ \tilde{\gamma} \circ \Phi ( x) }}_{g \cdot x}) \omega_{g \cdot x} ( dg_x v) \\
&= \iota( dg_x {X_{ \tilde{\gamma} \circ \Phi ( x) }}_{ x}) \omega_{g \cdot x} ( dg_x v) \\
&= \iota(  {X_{ \tilde{\gamma} \circ \Phi ( x) }}_{ x}) \omega_{ x} ( v) =  \iota(Y) \omega _x (v)\\
\end{split}
\end{equation*}

(2) is true because $Y$ is tangent to the fibres and the fibres are isotropic.

\end{proof}
\end{prop}
Since it's basic, we have $\iota(Y) \omega = \Phi^* \nu$ for some 1-form $\nu$ on $U_0$. Write $f_t = \exp(tY)$.

\begin{cor}
$f_t^*( \iota(Y) \omega ) = \iota(Y) \omega$
\begin{proof}
Since $M_{U_0}$ is dense in $M_U$, it suffices to check there. By the above, we have there that $f_t^* \iota(Y) \omega = f_t^* \Phi^* \nu = (\Phi \circ f_t)^* \nu = \Phi^* \nu = \iota(Y) \omega$. 

\end{proof}
\end{cor}

\begin{cor}
$\iota(Y) \omega$ is exact.
\begin{proof}
$ \frac{d}{dt} f_t^* \omega = f_t^* \mathcal{L}_Y \omega = f_t^* d \iota(Y) \omega = d f_t^* \iota(Y) \omega = d \iota(Y) \omega$. Integrating from 0 to 1, have $0 = f_1^* \omega - f_0^* \omega = d \iota(Y) \omega$, since $f_0 = id$ and $f_1 = f$ is a symplectomorphism.  Therefore, have shown that $\iota(Y) \omega$ is closed, so defines a cohomology class.  

To see exactness of $\iota(Y) \omega$, we need to show that its cohomology class is zero.  Recall that we chose $U$ such that $M_U$ is a deformation retract of $M_\alpha$, so the inclusion $j: M_{\alpha} \to M_U$ induces an isomorphism in cohomology $j^*: M_u \to M_{\alpha}$.   Since $Y$ is tangent to $M_{\alpha}$, and $M_{\alpha}$ is isotropic, $j^* (\iota(Y) \omega) = \iota(Y) j^* \omega = 0$, so $\iota(Y) \omega$ is exact.

\end{proof}
\end{cor}

\begin{cor}\label{ltsurj}
The map $\Lambda: C^{\infty} \to \mathcal{H}$ is surjective.
\begin{proof}
Using the notation of the above corollaries, $\iota(Y) \omega = dh$ for some $h \in C^{\infty}(M_U)$. Since $\iota(Y) \omega$ is G-invariant, can choose $h$ to be $G$-invariant (e.g. by averaging). Then $f$ is the time-1 Hamiltonian flow of the $G$-invariant function $h$.

\end{proof}
\end{cor}

\begin{proof}[Proof of Theorem \ref{delzant}]
By the above, we have a short exact sequence of sheaves of abelian groups  $0 \to \underline{\ell \times \mathbb{R}} \to C^{\infty} \to \mathcal{H} \to 0$, inducing a long exact sequence in cohomology.  $C^{\infty}$ is ``flabby'', so $H^i(\Delta, C^{\infty})=0, \forall i>0$.  $\Delta$ is contractible, so $H^i(\Delta,  \underline{\ell \times \mathbb{R}}) = 0, \forall i > 0$.  The long exact sequence then gives $H^1( \Delta, \mathcal{H}) = 0$, which completes the argument.

\end{proof}

\section{Canonical forms for near-symplectic 4-manifolds}

A \emph{near-symplectic structure} on a compact 4-manifold $M$ is a closed 2-form $\omega$ which is self-dual and harmonic with respect to some metric $\rho_{\omega}$, and is transverse to the zero section of the bundle $\wedge_2^{+,\rho}$ of self-dual 2-forms. By transversality, the zero set $Z_{\omega}$ of $\omega$ is  a 1-manifold, i.e. a disjoint union of circles $C_i$.  We call $Z_{\omega}$ the \emph{vanishing locus}. Since $\omega_p \wedge \omega_p = \omega_p \wedge *\omega_p = 0$ iff $\omega_p =0$, $\omega$ is symplectic on $M \setminus Z_{\omega}$ (the \emph{symplectic locus}).  It's a result due to Honda that if $b_2^+(M)>0$, for generic pairs $(\rho, \omega)$ with $\omega$ $\rho$-self-dual and harmonic, $\omega$ is transverse, i.e. $(M,\omega)$ is near-symplectic.

\begin{rem}
Auroux et. al. (\cite{Au}) give an equivalent definition of a near-symplectic structure that's independent of a Riemannian metric and show that a metric with respect to which the form is self-dual can always be constructed.  Their analysis is similar to that in our Appendix.
\end{rem}

In \cite{Ho}, Honda proves that near each component $C_i$ of $Z_{\omega}$, there is a neighborhood that is symplectomorphic to one of two canonical models $(S^1 \times D^3, \omega_A)$ and $(S^1 \times D^3, \omega_B)$. I'll present his argument in this chapter. The main tool in the proof is a Darboux-Moser type theorem for near-symplectic structures which will be extended to the equivariant case in Chapter 4.  The theorem relies on the existence of a canonical splitting of the normal bundle $NC$, which I'll describe first.

\subsection{Normal bundle splittings and standard forms}
Assume $M$ is oriented; then so is a neighborhood of $C$, $N(C)$, and so is the normal bundle to $C$, $NC$.  $\pi_1(BSO(3)) =0$, so this bundle is trivial and we can choose a $\rho_{\omega}$-orthonormal frame for $NC$ along $C$. Exponentiating the frame with respect to some metric ($\rho_{\omega}$ works, but so does any other one, e.g. a $G$-invariant one in the presence of a $G$-action), we obtain a diffeomorphism $S^1 \times  D^3 \to ^{\psi} N(C)$. Let $(\theta, x_1, x_2, x_3)$ be such coordinates on $S^1 \times D^3$.  Then in the chart $\psi$, the tangent vectors $\{ \frac{\partial}{ \partial \theta}, \frac{\partial}{\partial x_i} \}$ form an oriented $\rho_{\omega}$-orthonormal basis at all points $(\theta, 0)$.

Since $\omega (\theta, 0) = 0$, we can Taylor-expand $\omega$ in the coordinates $\psi$ to write
\begin{equation}
\begin{split}
\omega &= L_1(\theta,x)( d\theta dx_1 + dx_2 dx_3) \\
&+ L_2(\theta,x)( d\theta dx_2 + dx_3 dx_1) \\
&+ L_3(\theta,x)( d\theta dx_3 + dx_1 dx_2) + Q
\end{split}
\end{equation}
where $L_i(\theta,x) = \sum_{i=1}^3 L_{ij}(\theta)x_j$ are linear in $x$ and $Q$ is quadratic or higher in $x$. Note that this particular form holds since $\omega$ is $\rho_{\omega}$-self-dual with $\{ \frac{\partial}{ \partial \theta}, \frac{\partial}{\partial x_i} \}$ an oriented $\rho_{\omega}$-orthonormal basis at all points $(\theta, 0)$.

Using $d\omega=0$, calculating using the above expression, and equating 0th order terms, we obtain

\begin{subequations}
\begin{equation}
\frac{\partial L_1}{\partial x_2} - \frac{\partial L_2}{\partial x_1} = 0, \frac{\partial L_1}{\partial x_3} - \frac{\partial L_3}{\partial x_1} = 0, \frac{\partial L_2}{\partial x_3} - \frac{\partial L_3}{\partial x_2} = 0 \\
\end{equation}
\begin{equation}
\frac{\partial L_1}{\partial x_1} +\frac{\partial L_2}{\partial x_2} +  \frac{\partial L_3}{\partial x_3} = 0,  
\end{equation}
\end{subequations}

Since $\frac{\partial L_i}{\partial x_j} = L_{ij}(\theta)$, this shows that the matrix $\{L_{ij}\}_{i,j=1}^3$ is traceless and symmetric.  Thus $L_{ij}$ is diagonalizable.  

\begin{prop}
For $\omega$ transverse to the zero section of $\wedge^+_{\rho_{\omega}}$, $L_{ij}$ has full rank.
\begin{proof}
The fibre of  $\wedge^+_{\rho_{\omega}}$ at a point $x \in M$ has dimension 3, and since $\omega=0$ along $C$, the image of $\partial \omega |_{N_xC}$ must span the fibre.  $L_{ij}$ is this derivative.

\end{proof}
\end{prop}

By the Proposition, $L_{ij}$ has no zero eigenvalues. Thus it must have two positive eigenvalues and one negative eigenvalue (or vice versa).  

I will now describe these eigenspaces in an invariant way as subspaces of the normal bundle $NC$ using only the metric $\rho_{\omega}$.

Let $A=\sum a_i \frac{\partial}{\partial x_i}, B = \sum b_i \frac{ \partial }{ \partial x_i } \in N_{(\theta, 0)}C$  be two tangent vectors in the above coordinates. Then
\begin{equation}\label{Lform}
\begin{split}
\sum_{i,j} a_j b_i L_{ij} &= \sum_i b_i( \sum_j a_j \iota(\frac{\partial}{\partial x_j} ) d[\iota(\frac{\partial}{\partial x_k})\iota(\frac{\partial}{\partial x_{\ell}})\omega])\\
&= \iota(A) \sum_i b_i d[\iota(\frac{\partial}{\partial x_k})\iota(\frac{\partial}{\partial x_{\ell}})\omega]\\
\end{split}
\end{equation}
where $\{ \frac{\partial}{\partial x_i} , \frac{\partial}{\partial x_k}, \frac{\partial}{\partial x_{\ell}} \}$ is an oriented basis of $N_{(\theta, 0)}C$.

\begin{prop}
Let $V \cong \mathbb{R}^3$ be a vector space with the standard inner product and orientation, and let $\omega$ be a skew-symmetric bilinear form on $V$.  Then the map $ v \in V \mapsto ^q \omega ( v', v'')$, where $\{v, v', v''\}$ is an oriented orthogonal basis of $V$ and $|v'| = |v''| = |v|^{1/2}$, is well-defined and linear.
\begin{proof}
$q$ is well-defined because any orientation-preserving orthogonal transformation $(v', v'') \mapsto (\tilde{v}', \tilde{v}'')$ leaves $\omega$ invariant. One can check that if $(x,y,z)$ are standard coordinates on $\mathbb{R}^3$, $\omega = a  dx\wedge dy + b dy \wedge dz + c dx \wedge dz$, and $v = (x_0, y_0, z_0)$, then $q(v) = az_0 + bx_0 - cy_0$.

\end{proof}
\end{prop}

\begin{cor}\label{splitting}
Given $\rho_{\omega}$, there is a natural orthogonal splitting of the normal bundle $NC$ into a 2-dimensional subbundle and a line bundle.
\begin{proof}
By the previous proposition we can define a bilinear form $H: N_pC \times N_pC \to \mathbb{R}$ by 
\begin{equation}
H(v,w) = \iota(v) (d (q(\tilde{w}))_p
\end{equation}
where $q$ is as in the proposition and $\tilde{w}$ is any extension of $w$ to a vector field near $p$.  (H is independent of the choice of $\tilde{w}$ by vanishing of $\omega$ at $C$.) By \eqref{Lform}, the associated map $\tilde{H}: N_pC \to N_pC$ obtained via $\rho_{\omega}$ is represented by the matrix $L_{ij}$, and so induces an orthogonal splitting of the normal bundle $NC$ into a 2-dimensional subbundle and a line bundle which are the spans of the positive and negative eigenspaces.

\end{proof}
\end{cor}
 
 Since $NC \cong S^1 \times D^3$ is trivializable, we can classify such splittings by maps $S^1 \to \mathbb{RP}^2$. Up to homotopy, these are classified by $\pi_1( \mathbb{RP}^2) = \mathbb{Z}/ 2 \mathbb{Z}$. We can distinguish these splittings by whether the line bundle is orientable or not.  
 
 \begin{prop}[\cite{Ho}]
There are model near-symplectic structures on $S^1 \times D^3$ with vanishing locus  $x=0$ and self-dual with respect to the flat metric that induce both types of splittings. 
 \begin{proof}
 Representatives $\omega_A$ and $\omega_B$ are defined as follows.
 
The oriented splitting: On $S^1 \times D^3$, set 
\begin{equation}
\begin{split}
\omega_A &= x_1( d\theta dx_1 + dx_2 dx_3) \\
&+ x_2(d \theta d x_2 + d x_3 d x_1)\\
&- 2x_3(d \theta d x_3 + d x_1 d x_2)
\end{split}
\end{equation}
 here $L_{ij}(\theta) = \text{diag}(1,1,-2)$ with fixed positive and negative eigenspaces.
 
 The unoriented splitting: Set $\Omega = \omega_A$ on  $[0,2\pi] \times D^3$.  Then glue $\{2 \pi\} \times D^3 \to^{\phi} \{0\} \times D^3$ by $ \theta \mapsto \theta - 2 \pi$, $x_1 \mapsto x_1$, $x_2 \mapsto -x_2$, $x_3 \mapsto -x_3$.  Then $\phi^* \Omega = \Omega$ so $\Omega$ induces a form $\omega_B$ on the quotient.
 
 \end{proof}
 \end{prop}
 
 \subsection{Contact boundaries and Reeb flow}
 
 \begin{prop}[\cite{Ho}] 
Both models $(S^1 \times D^3, \omega_A)$ and $(S^1 \times D^3, \omega_B)$ admit compatible contact structures on their boundaries $S^1 \times S^2$.
\begin{proof}
For (A), consider the 1-form
\begin{equation}
\lambda = -\frac{1}{2}(x_1^2 + x_2^2 - 2x_3^2) d \theta + x_2 x_3 dx_1 - x_1 x_3 dx_2
\end{equation}
We have $d \lambda = \omega_A$.  Let $i: S^1 \times S^2 \to S^1 \times D^3$ be the inclusion.  Then since $\sum_i x_i dx_i = 0$ on $T(S^1 \times S^2)$, $i^*(\lambda \wedge d \lambda) \neq 0$ iff $\lambda \wedge d\lambda \wedge \sum_i x_i dx_i \neq 0$ near $S^1 \times S^2$.  But
\begin{equation}
\lambda \wedge d\lambda \wedge \sum_i x_i dx_i = -(\frac{1}{2}(x_1^2 + x_2^2)(x_1^2 + x_2^2 +2x_3^2) + 2x_3^4) d \theta dx_1 dx_2 dx_3
\end{equation}
so is nonzero where required.  

For case (B), the proof is the same after gluing.
\end{proof}
\end{prop}
 
 To conclude this section, I reproduce Honda's description of the Reeb vector fields on $N(C)$.   Note that the flat metric on $S^1 \times D^3$ is given by $\rho(x,y) =  \frac{1}{c}\omega_A (x, Jy)$, where $J = -\frac{1}{c} A$, $A$ is the matrix representation of $\omega_A$, i.e.
 \begin{equation}
A = \begin{pmatrix}
0 &x_1 & x_2 & -2x_3 \\
-x_1 & 0 & -2x_3 & -x_2 \\
-x_2 & 2x_3 & 0 & x_1 \\
2x_3 & x_2 & -x_1 & 0
\end{pmatrix} 
 \end{equation}
 and $c = \sqrt{x_1^2 + x_2^2 +4x_3^2}$.  By compatibility, requiring the Reeb vector field $X$ to be in $\text{ker}(i^* d\lambda)$ is equivalent to requiring it to be in the image under $J$ of the $\rho$-normal to $S^1 \times S^2$, i.e. up to scalars, 
 \begin{equation}
 X = J(\sum_i x_i \frac{\partial}{\partial x_i}) =  \frac{-1}{\sqrt{x_1^2 + x_2^2 +4x_3^2}}( (x_1^2 +x_2^2 -2x_3^2) \frac{\partial}{\partial \theta} - 3 x_2 x_3 \frac{\partial}{\partial x_1} + 3 x_1 x_3 \frac{\partial}{\partial x_2})
 \end{equation}
 Normalizing by $\lambda(X) =1$ gives
 
 \begin{equation}
 X =  \frac{1}{f}( (x_1^2 +x_2^2 -2x_3^2) \frac{\partial}{\partial \theta} - 3 x_2 x_3 \frac{\partial}{\partial x_1} + 3 x_1 x_3 \frac{\partial}{\partial x_2})
 \end{equation}
 where $f = -\frac{1}{2}[(x_1^2 + x_2^2)(x_1^2 + x_2^2 +2x_3^2) + 4x_3^4] $.

Setting $r^2 = x_1^2 + x_2^2, \beta = \arctan(x_2/x_1)$, we can rewrite the vector field as 

\begin{equation}
 X =  \frac{1}{f}( (r^2 -2x_3^2) \frac{\partial}{\partial \theta} + 3 x_3  \frac{\partial}{\partial \beta} )
 \end{equation}
 where $f = -\frac{1}{2}[(r^2)(r^2 +2x_3^2) + 4x_3^4] $.

Thus the flow preserves the $x_3$ coordinate and  $r^2$, and  rotates in the $(x_1, x_2)$-plane and along $S^1$. On $S^1 \times S^2$, we have $r^2 + x_3^2 = k$ (usually $k=1$, but can take any $k \neq 0$), so the flow can be written in the form
\begin{equation}
\begin{split}
x_1(t)& = \sqrt{k-x_3^2} \cos( R_1(x_3) t) \\
x_2(t)& = \sqrt{k-x_3^2} \sin( R_1(x_3) t) \\
x_3(t)& = x_3(0) \\
\theta(t)& = R_2(x_3)t + c
\end{split}
\end{equation}
where $R_i$ are functions of $x_3$. Specifically, since $-2 f = k^2 +3x_3^4$, 
\begin{equation}
\begin{split}
R_1 &= -2 \frac{3 x_3 }{ k^2 + 3x_3^4} \\
R_2 &= -2 \frac{r^2 - 2x_3^2 }{k^2 + 3x_3^4}
\end{split}
\end{equation}

We can now consider the closed orbits of the Reeb flow. Note that for $x_3 = 0$, $R_1 = 0$, and the closed orbit is of the form $(x_1,x_2, 0) = \text{constant}$, i.e. flow along the $\theta$ direction.   Similarly, for $r=0$, the flow is along the $\theta$ direction.  For $r^2 -2x_3^2=0$, $R_2=0$, so the closed orbit is of the form $(x_3, r, \theta) = \text{constant}$, i.e. flow along the $\beta$ direction. The other closed orbits occur when $R_1 / R_2 \in \mathbb{Q}$. 
\begin{rem}
Note that the cases $x_3 =0$ and $r=0$ correspond respectively to the stable and unstable gradient directions in the Morse-Bott theory for the function $x_1^2 + x_2^2 - 2 x_3^2 = r^2-2x_3^2$, which is the moment for the $S^1$ action given by rotation in $\theta$ and the numerator in $R_2$.   The significance of the numerator in $R_1$ is unclear.
\end{rem}

\begin{rem}
Honda (\cite{Ho}) also proves that the contact structures induced on the boundaries $S^1 \times S^2$ by the contact forms $\lambda_A$ and $\lambda_B$ are both overtwisted and distinct, but I'll skip the proof because I won't make use of it later.
 \end{rem}
 
 \subsection{Honda-Moser theorems}
Let $\{\omega_t\}, t \in [0,1]$ be a smooth family of self-dual harmonic 2-forms with respect to metrics $\rho_t$, transverse to the zero sections of their respective bundles of self-dual forms, such that the number of components of $Z_{\omega_t}$ is constant so that we can identify all $Z_{\omega_t}$ via isotopy. For simplicity, assume $Z_\omega = C$ is constant.  Assume further that (i) $[\omega_t] \in H^2(M; \mathbb{R})$ is constant, and (ii) $[\omega_t] \in H^2(M,C;\mathbb{R})$ is constant.  

\begin{thm}[Global Honda-Moser \cite{Ho}]\label{hondaglob}
Under the above assumptions, there exists a 1-parameter family $f_t$ of $C^0$-homeomorphisms of $M$, smooth away from $C$ and fixing $C$, such that $f_t^* \omega_t = \omega_0$.
\end{thm}

As in the proof of the Moser-Darboux theorem, the requirement $f_t^* \omega_t = \omega_0$, $f_0 = \text{Id}$, implies that $f_t^*(d \iota(X_t) \omega_t) + f_t^*( \frac{d \omega_t}{dt}) = 0$, where $X_t = \frac{df_t}{dt}$.  The proof thus reduces to choosing a 1-form $\eta_t$ satisfying $d \eta= \frac{d \omega_t}{dt} $ such that the equation $\iota(X_t) \omega_t = - \eta_t$ defines a vector field $X_t$ which is sufficiently continuous and zero along $Z_{\omega}$.  The complication in the near-symplectic case is that $\omega_t$ is degenerate along $Z_{\omega}$ so can't be smoothly inverted.  To deal with this complication, we will choose $\eta$ very carefully.

\begin{lem}\label{hondaexact}
There exists a smooth family of 1-forms $\tilde{\eta}_t$ such that $\frac{d \omega_t}{dt} = d \tilde{\eta}$ and $i^* \tilde{\eta}_t$ is exact, where $i: C \to M$ is the inclusion.
\begin{proof}
Recall the relative cohomology exact sequence:
\begin{equation*}
H^1(M;\mathbb{R}) \to ^{i^*} H^1(C;\mathbb{R}) \to ^{\delta} H^2(M,C;\mathbb{R}) \to H^2(M;\mathbb{R})
\end{equation*}
In deRham cohomology, for $[\alpha] \in H^1(C, \mathbb{R})$, $\delta [\alpha] = [d \tilde{\alpha}]$ where $\tilde{\alpha}$ is any extension of $\alpha$ to a 1-form on $M$.

By assumption (i), there exists a smooth family of 1-forms $\hat{\eta}_t$ such that $\frac{d \omega_t}{dt} = d \hat{\eta}_t$. $\hat{\eta}_t$ is clearly an extension of $i^* \hat{\eta}_t$ to $M$, and by condition (ii), $[d \hat{\eta}_t] = 0 \in H^2(M,C;\mathbb{R})$.  Therefore, by exactness, there exists $[\beta_t] \in H^1(M;\mathbb{R})$ such that $[i^* \hat{\eta}_t] = [i^* \beta_t]$. Setting $\tilde{\eta}_t = \hat{\eta}_t - \beta_t$ gives the required 1-form.

\end{proof}
\end{lem}

\begin{lem}\label{hondasecond}
There exists a smooth family of 1-forms $\eta_t$ such that $\frac{d \omega_t}{dt} = d \eta_t$ and $\eta_t = 0$ at $Z_{\omega}$ \emph{up to second order}.
\begin{proof}

Let $\tilde{\eta}_t$ be as in the previous lemma and let $f_t$ be a smooth family of functions such that $i^* \tilde{\eta}_t = df_t$. Choose coordinates on $N(C) \cong S^1 \times D^3$ via exponentiating a  $\rho_{\omega}$-orthonormal frame along $C$ as above.  Set
\begin{equation*}
f_t( \theta, x_1, x_2, x_3) = f_t (\theta) + \sum_i \tilde{\eta}_i(\theta,0) x_i + \frac{1}{2} \sum_{i,j} \frac{\partial \tilde{\eta}_i}{\partial x_j}(\theta,0)x_i x_j
\end{equation*}
where $\tilde{\eta} = \tilde{\eta}_{\theta}d\theta + \sum_i \tilde{\eta}_i dx_i$ on N(C). Then
\begin{equation*}
\begin{split}
df_t(\theta, x_1, x_2, x_3) &= \frac{\partial f_t}{\partial \theta}(\theta) d \theta + \sum_i \frac{ \partial \tilde{\eta}_i}{\partial \theta}(\theta, 0) x_i d \theta \\
&+ \sum_i \tilde{\eta}_i(\theta, 0) d x_i + \frac{1}{2} \sum_{i,j} \frac{\partial \tilde{\eta}_i}{\partial x_j}(\theta, 0)(x_i dx_j + x_j dx_i) \\
&+ O(x^2)
\end{split}
\end{equation*}
Notice that 
\begin{equation*}
\begin{aligned}
\frac{\partial f}{\partial \theta}(\theta) =& \tilde{\eta}_{\theta}(\theta,0) \\
d \tilde{\eta}_t (\theta, 0) =& \frac{d\omega_t}{dt}(\theta,0) =0 \\
\implies & \frac{\partial \tilde{\eta}_{\theta}}{\partial x_i}(\theta, 0) = \frac{\partial \tilde{\eta}_i}{\partial \theta}(\theta, 0) \\
&  \frac{\partial \tilde{\eta}_i}{\partial x_j}(\theta, 0) = \frac{\partial \tilde{\eta}_j}{\partial x_i}(\theta, 0) \\
\end{aligned}
\end{equation*}
Substituting, we obtain
\begin{equation}
\begin{split}
df_t(\theta, x_1, x_2, x_3) &= (\tilde{\eta}_{\theta}(\theta,0) + \sum_i \frac{\partial \tilde{\eta}_{\theta}}{\partial x_i} (\theta,0) x_i) d\theta \\
&+ \sum_i(\tilde{\eta}_i(\theta,0) + \sum_j(\frac{\partial \tilde{\eta}_i}{\partial x_j}(\theta,0)x_j))dx_i \\
&+ O(x^2)
\end{split}
\end{equation}
Setting $\eta_t = \tilde{\eta}_t - df_t$ (after damping $f_t$ to zero away from C by a cutoff function, so that it extends by zero to all of $M$) we obtain the required 1-form.

\end{proof}
\end{lem}

\begin{proof}[Proof of theorem \ref{hondaglob} ]
Define the vector field $X_t$ by $\iota(X_t) \omega = - \eta_t$, where $\eta_t$ is as in the previous lemma. As a matrix, we've seen that $\omega_t$ corresponds to
\begin{equation*}
A(\theta,x) = 
\begin{pmatrix}
0 & L_1 & L_2 & L_3 \\
-L_1 & 0 & L_3 & -L_2 \\
-L_2 & -L_3 & 0 & L_1 \\
-L_3 & L_2 & -L_1 & 0
\end{pmatrix} (\theta,x)
+ Q(\theta,x)
\end{equation*}
where $Q$ is quadratic or higher in $x$, and this is in the coordinates coming from the $\rho_{\omega}$-orthonormal frame along $C$. Setting $X_t = a_t \frac{\partial}{\partial \theta} + \sum_i a_i \frac{\partial}{\partial x_i}$, we have $(a_\theta, a_1, a_2, a_3) A = - (\eta_{\theta}, \eta_1, \eta_2, \eta_3)$, or
\begin{equation*} 
\begin{split}
(a_\theta, a_1, a_2, a_3) &= - (\eta_{\theta}, \eta_1, \eta_2, \eta_3) A^{-1} \\
&= (\eta_{\theta}, \eta_1, \eta_2, \eta_3) \frac{1}{L_1^2 + L_2^2 + L_3^2} \begin{pmatrix}
0 & L_1 & L_2 & L_3 \\
-L_1 & 0 & L_3 & -L_2 \\
-L_2 & -L_3 & 0 & L_1 \\
-L_3 & L_2 & -L_1 & 0
\end{pmatrix}  + Q'
\end{split}
\end{equation*} 
where $Q'$ is second order or higher in $x$.

Notice that by nondegeneracy of $L_{ij}$, $L_1^2 + L_2^2 + L_3^2 \neq 0$ unless $x=0$.  Given this nondegeneracy, the expression above has leading term of order 1 in $x$. Thus $|X_t| < k|x|$ near $C$, $X_t$ is smooth elsewhere, and the flow fixes $C$.

\end{proof} 

\begin{thm}[Local Honda-Moser \cite{Ho}]\label{hondaloc}
Let $(M,\omega)$ be a near-symplectic manifold. Then near each component $C_i$ of $Z_{\omega}$, there is a neighborhood that is symplectomorphic to one of the two models $(S^1 \times D^3, \pm \omega_A)$ or $(S^1 \times D^3, \pm \omega_B)$.
\begin{proof}
Assume that $\omega$ is such that $NC$ splits in the oriented manner. The proof for the unoriented case is similar (working in the nonreduced space of the unoriented model). We can choose coordinates via an $\rho_{\omega}$-orthonormal frame such that  $\{ \frac{\partial}{\partial x_i} \}_{i=1,2}$ span the two dimensional subbundle along $C$ and $\{ \frac{\partial}{\partial x_3} \}$ spans the one-dimensional subbundle.  Then in these coordinates,
\begin{equation*}
\begin{split}
\omega &= (L_{11}(\theta)x_1 + L_{12}(\theta)x_2)( d\theta dx_1 + dx_2 dx_3) \\
 &+ (L_{21}(\theta)x_1 + L_{22}(\theta)x_2)( d\theta dx_2 + dx_3 dx_1) \\
  &+ L_{33}(\theta)x_3( d\theta dx_3 + dx_1 dx_2) + Q 
\end{split}
\end{equation*}
where $(L_{ij})_{i,j = 1,2}$ is positive-definite and $L_{33} < 0$ (for the opposite case, change signs in $\omega_A$ below to match).  Using these $\omega$-adapted coordinates, write 
\begin{equation*}
\begin{split}
\omega_A &= x_1(d\theta dx_1 + dx_2 dx_3) \\
&+ x_2(d\theta dx_2 + dx_3 dx_1) \\
&- 2x_3(d\theta dx_3 + dx_1 dx_2) \\
\end{split}
\end{equation*}
Now, set $\omega_t = t\omega + (1-t) \omega_A$.  We can still define $L_{ij}(t)$ as before, only we do so in this fixed coordinate system. Since we defined $\omega_A$ so that its eigenspaces correspond to those of $\omega$, $L_{ij}(t)$ is nondegenerate for all $t \in [0,1]$ and has the same form as in the global theorem. Finally, on a tubular neighborhood $N(C)$, the cohomological conditions of the global theorem are trivially satisfied, so the proof of the global theorem carries through otherwise unmodified.

\end{proof}
\end{thm}

\section{Near-symplectic toric 4-manifolds}

I now consider the case in which $T^2$ acts effectively on a near-symplectic 4-manifold.  This chapter is inspired by the work of Gay-Symington ($\cite{Gay-S}$) but takes a different approach.  In particular, I make the simplifying assumption throughout that there is a global Hamiltonian $T^2$ action rather than only a locally toric structure. I use the existence of a metric as in Condition \ref{invariantmetric} to show that near each component of $Z_{\omega}$ the toric structure is of a certain standard form.

\begin{defn}[\cite{Gay-S}]
A smooth $T^2$ action on a near-symplectic manifold $(M,\omega)$ is \emph{Hamiltonian} if there exists a smooth map $\Phi: M \to \mathfrak{t}^*$ such that $\Phi$ is a moment map for the action on $M \setminus Z_{\omega}$.  A near symplectic $4$-manifold is \emph{toric} if it has an effective Hamiltonian $T^2$-action.
\end{defn}

From this definition, it's clear that the local analysis of the moment map in Chapters 1 and 2, including the local canonical forms, local convexity, and Delzant conditions, carry through unchanged on $M \setminus Z_{\omega}$.

\subsection{Equivariant Honda-Moser theorems}

To obtain a canonical local form for the $T^2$ action near $Z_{\omega}$, it's convenient to extend the results of Chapter 3 to be equivariant.  It is easy to extend Lemmas \ref{hondaexact}, \ref{hondasecond} and Theorem \ref{hondaglob} to the equivariant setting, as follows.  

\begin{lem}
In Lemmas \ref{hondaexact}, \ref{hondasecond} and Theorem \ref{hondaglob} we can choose the one-forms $\tilde{\eta}$ and $\eta$ to be $G$-invariant and the diffeomorphism $f$ to be $G$-equivariant, for any connected compact Lie group $G$ that acts symplectically on $(M, \omega_t)$.
\begin{proof}
Since averaging over $G$ preserves cohomology classes and $\omega_t$ is invariant,  $\tilde{\eta}$ may be averaged at the end of the proof of Lemma \ref{hondaexact} to obtain an invariant form satisfying all the conditions. Similarly, since $\omega_t$ is invariant and $G$ is compact and fixes $Z_{\omega}$, $\eta$ can be averaged at the end of the proof of Lemma \ref{hondasecond} to obtain an invariant form as required, preserving the estimate on its vanishing.  Finally, since $\eta_t$ and $\omega_t$ are $G$-invariant, so is the vector field $X_t$ in the proof of Theorem \ref{hondaglob}, and its flow $f$ is $G$-equivariant.
\end{proof}
\end{lem}

It is less easy to extend Theorem \ref{hondaloc} in a useful way, making use of a $T^2$ action.  For this purpose I state the following condition explicitly: 

\begin{cond}\label{invariantmetric}
Assume there exists a Riemannian metric $\rho$ such that $\omega$ is self-dual and transverse with respect to $\rho$ \emph{and $\rho$ is $T^2$-invariant}.
\end{cond}

\begin{rem}
This condition is always satisfied for any near-symplectic $\omega$.  A proof of this fact due to D. Auroux is explained in the Appendix.
\end{rem}

The canonical form $\omega_A$ on $S^1 \times D^3$ can be given a Hamiltonian $T^2$ action which satisfies the above condition with the flat metric:

\begin{exmp}[The ``standard fold'']\label{canfoldmodel}
Let $\omega_A = x(d\alpha dx + dy dz) + y(d\alpha dy + dz dx) -2z(d\alpha dz + dx dy)$ on $S^1 \times D^3 \cong (\alpha, x, y ,z)$. Setting $\theta = \arctan(y/x), r^2 = x^2 + y^2$, we have $\omega_A = -2z(d \alpha dz + r dr d\theta) - r dr d\alpha + r^2 d\theta dz$. Then\\
1. $\omega_A$ is invariant under the $T^2$ action $(t_1, t_2) \cdot (\alpha, r, \theta, z) = (\alpha +t_1, r, \theta + t_2, z)$\\
2. A moment map for the action is $(\alpha, r, \theta, z) \mapsto ^{\Phi_0} (z^2 - \frac{1}{2}r^2, zr^2)$ \\
3. $\omega_A$ is self-dual and transverse with respect to the flat metric.\\
4. The $T^2$ action preserves the flat metric.
\begin{proof}[Proof of 2]
Setting $p_1 = z^2 - \frac{1}{2}r^2, q_1=\alpha, p_2 = zr^2, q_2 = \theta$, we calculate $\omega_A = dp_1 \wedge dq_1 + dp_2 \wedge dq_2$, as desired.
\end{proof}
\end{exmp}

It is worth pointing out several features of this model in detail.  First, the orbit space $B = M / T^2$ can be identified with the half-plane $H = \{ (x_1, x_2) \in \mathbb{R}^2 | x_2 \geq 0 \}$ via the map $(\alpha, r, \theta, z) \mapsto (z, r^2)$.  Under this identification, $\partial B =  \{ (x_1, x_2) \in H | x_2 = 0 \}$.  Let $\pi: M \to B$ be the quotient map, and let $S^1_1 \times S^1_2 \cong T^2$ be the standard splitting of $T^2$ into two circle subgroups. Since $r=0$ there, each point in $\pi^{-1} ( \partial B)$ has stabilizer $S^1_2$, and lies on a circle orbit generated by $S^1_1$.  On $\text{int}(B)$, each point has preimage a full $T^2$ orbit with trivial stabilizer.

Let $\mathfrak{t}_1 \oplus \mathfrak{t}_2$ be the splitting of the Lie algebra $\mathfrak{t}$ induced by the group splitting $T^2 \cong S^1_1 \times S^1_2$. It induces a splitting $\mathfrak{t}^* \cong \mathfrak{t}^*_1 \oplus \mathfrak{t}^*_2$.  The moment map $\Phi_0$ descends to a map $\phi_0: B \to \mathfrak{t}^*$, given in the coordinates $(x_1, x_2)$ on $H$ and  $\mathfrak{t}^*_1 \oplus \mathfrak{t}^*_2$ on $\mathfrak{t}^*$ by $(x_1, x_2) \mapsto ^{\phi_0} (x_1^2 - \frac{1}{2} x_2, x_1 x_2)$.  This map is, in the terms of Gay-Symington (\cite{Gay-S}), a \emph{fold}, that is, it satisfies the following properties:\\
0.  $\phi_0: H \to \mathbb{R}^2$ is smooth. \\
1. $\phi_0(0,0) = (0,0)$. \\
2. ${\phi_0}_{H \setminus \{(0,0)\}}$ is an immersion. \\
3. $\phi_0$ maps both $\{(x_1,0)|x_1>0\}$ and $\{(x_2,0)|x_2<0\}$ diffeomorphically onto $\{(p_1,0)|p_1>0\}$. \\
4. $\phi_0$ maps $\{(x_1,x_2)|x_2>0\}$ diffeomorphically onto $\mathbb{R}^2 \setminus \{(p_1,0)|p_1>0\}$.

\begin{figure}
\PSbox{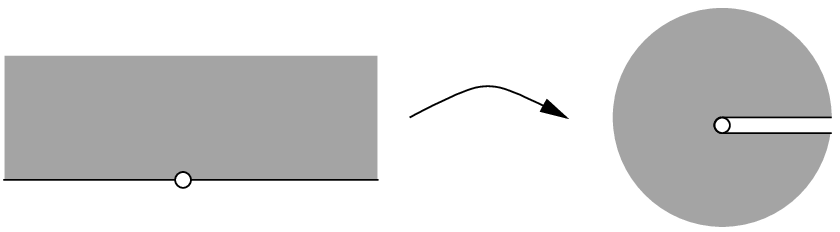}{3.4 in}{1.0 in}
\caption{A fold.}
\label{foldfig}
\end{figure}

A familiar example of a fold is the complex map $z \mapsto z^2$, restricted to $H \subset \mathbb{C}$. A fold is illustrated in Figure \ref{foldfig}, with the double $p_1$-axis drawn as two parallel lines for the purpose of illustration.

Here $(p_1, p_2)$ are coordinates on $\mathfrak{t}^*$.  The above description of a fold means that $\phi_0$ maps $\partial B$ to the positive $p_1$-axis. Along the positive $p_1$-axis, each point $x$ has two pre-images under $\phi_0$, which lie in $\partial B$, i.e. two orbits map to $x$ and they are both circle orbits.  The origin $(p_1,p_2)=(0,0)$ has a unique pre-image, also a circle orbit in $\partial B$. Notice that the double image along the positive $p_1$-axis lies in $\mathfrak{t}^*_1 \times \{0\}$, i.e. in $\mathfrak{t}_2^0 \subset \mathfrak{t}^*$, the annihilator of $\mathfrak{t}_2$, the infinitesimal generator of the the stabilizer group of $\pi^{-1}(\partial B)$. 

Finally, for  a point not on the positive $p_1$-axis, the pre-image under $\phi_0$ is a unique point in $\text{int}(B)$ corresponding to a free $T^2$ orbit, and the set $\Phi_0^{-1}( \mathfrak{t}^* \setminus \{ (p_1, 0) | p_1 \geq 0 \})$ is equivariantly symplectomorphic to the manifold $(\mathfrak{t}^* \setminus \{ (p_1, 0) | p_1 \geq 0 \}) \times T^2$ with symplectic form $dp_1 \wedge dq_1 + dp_2 \wedge dq_2$, where $q_1$ is a coordinate on $S^1_1$ and $q_2$ is a coordinate on $S^1_2$, as shown in the example.  The existence of these standard coordinates away from the fold is what permits the easy construction of near-symplectic toric $4$-manifolds by patching together local models as in Sections 4.3 and 4.4.

\begin{rem}[Contact structure and Reeb flows]
We remark that the contact form $\lambda_A$ defined in Chapter 3 can be written  in the above notation as 
\begin{equation}
\lambda_A = -\frac{1}{2}(r^2 - 2z^2) d \alpha - z r^2 d \theta
\end{equation}
so it's $T^2$ invariant, and its coefficients are the components of the moment map $\Phi_0$. Indeed, this is just saying $\lambda_A = p_1 dq_1 + p_2 dq_2$, where $\omega_A = d\lambda_A = dp_1 \wedge dq_1 + dp_2 \wedge dq_2$ as above.  Similarly, the Reeb vector field $X$ can be written as
\begin{equation}
 X =  \frac{1}{f}( (r^2 -2z^2) \frac{\partial}{\partial \alpha} + 3 z  \frac{\partial}{\partial \theta} )
 \end{equation}
 where $f = -\frac{1}{2}[(r^2)(r^2 +2z^2) + 4z^4] $, so it's $T^2$ invariant and tangent to the fibres, and hence moment-preserving.  
  \end{rem}

As we will see in Theorem \ref{canfoldthm}, a consequence of Condition \ref{invariantmetric} is that near any component $C$ of $Z_{\omega}$, the manifold is equivariantly symplectomorphic to the above ``standard fold''.

\begin{prop}
Assuming Condition \ref{invariantmetric}, the $\rho$-normal bundle $NC$  and the splitting from Corollary \ref{splitting} are $T^2$ invariant.
\begin{proof}
The invariance of $NC$ follows since $T^2$ preserves $C$ and $\rho$ is $T^2$-invariant.\\
To see that the splitting is preserved, we want $H_{g \cdot p}(dg \cdot v_p, dg \cdot w_p) = H_p(v_p, w_p)$, for all $p \in C, g \in T^2$, and $H$ defined as in Corollary \ref{splitting}.
\begin{equation*}
H( dg \cdot v_p, dg \cdot w_p) = \iota(dg \cdot v_p) d q( \widetilde{ dg \cdot w_p }) 
\end{equation*}
Now, using the notation from Corollary \ref{splitting}, 
\begin{equation*}
\begin{split}
d(q(\widetilde{dg \cdot w})) & = d [ \omega( \widetilde{(dg \cdot w)'}, \widetilde{(dg \cdot w)''})] \\
&= d[ \omega( dg \cdot \widetilde{w'}, dg \cdot \widetilde{w''})]
\end{split}
\end{equation*}
because $T^2$ acts orthogonally and the quantity is independent of the extensions of the vector fields. So
\begin{equation*}
\begin{split}
\iota(dg \cdot v) d(q(\widetilde{dg \cdot w})) &= \iota(dg \cdot v) d[ \omega( dg \cdot \widetilde{w'}, dg \cdot \widetilde{w''})]_{g \cdot p} \\
&= \iota(v) d[ \omega( \widetilde{w'}, \widetilde{w''})]_p  = H(v,w)
\end{split}
\end{equation*}
where the second equality is by $T^2$ invariance of $\omega$ and the chain rule.
\end{proof}
\end{prop}

We are now ready to extend Theorem \ref{hondaloc} to the case with a $T^2$ action.

\begin{thm}\label{canfoldthm}
Assuming Condition \ref{invariantmetric}, any component $C$ of $Z_{\omega}$ has a neighborhood that's equivariantly symplectomorphic, up to an integral reparametrization of $T^2$, to the model in Example \ref{canfoldmodel}.
\begin{proof}

Fix a metric $\rho$ as in Condition \ref{invariantmetric}.  Choose $x \in C$. Consider $G_x$, the stabilizer group of $x$.  We first show that $G_x$ is a circle subgroup of $G = T^2$ and that the splitting of the normal bundle must be the oriented one.  Note that because $T^2$ acts symplectically, it must map $C$ to $C$.  

If $G_x = G$, then $G$ acts linearly on $T_xM$. Since $\rho$ is $G$-invariant, $G$ preserves the splitting of $NC$  and maps $C$ to $C$, $G$ acts on $T_xM$ as a subgroup of $O(1) \times O(1) \times O(2)$. By connectedness it must act as a subgroup of $SO(2)$, but this violates effectiveness, as there is no faithful representation $T^2 \to SO(2)$.

Since $G_x \neq G$ and is closed, it must have dimension $0$ or $1$.  If it had dimension $0$, the orbit $G/G_x$ would have dimension 2, but $G/G_x \subset C$, which is 1-dimensional.  Thus $G_x$ has dimension 1 and $C$ is the orbit $G \cdot x$.

Choose an orthonormal basis for $\{V_1, V_2, V_3\}$ of $N_xC$ such that $\{V_1,V_2\}$ span the 2-dimensional sub-space of the splitting of $N_xC$ and $V_3$ spans the 1-dimensional subspace.  Exponentiating these vectors by $\rho$, we obtain a slice for the $G$-orbit $C$ at $x$.  Since $G$ preserves the splitting of $NC$, $G_x$ acts on $N_xC$ as a subgroup of $O(1) \times O(2)$, and by effectiveness and the equivariant slice theorem, this representation is faithful, i.e. $G_x \subset O(1) \times O(2)$ as a subgroup.  By connectedness of $T^2$ and disconnectedness of $O(1)$, $G_x$ acts nontrivially on the 1-dimensional subspace iff the line bundle is nonorientable, i.e. we are in the unoriented splitting.  However, we must also have that $G_x$ preserves the orientation of $N_xC$ since it preserves the orientation of $C$ and acts symplectically.  The only 1-dimensional abelian subgroup of $(O(1) \times O(2)) \cap SO(3)$ is $ 1 \times SO(2)$, so $G_x$ is a circle subgroup of $T^2$ and the splitting must be the oriented splitting.

For any closed circle subgroup $G_x$ of $T^2$, we can choose a complement circle subgroup $H$ such that $T^2$ splits as $T^2 = G_x \times H$.  A consequence of effectiveness is that in these coordinates, the generators of the Lie algebras $\{\mathfrak{g}_x, \mathfrak{h}\}$ correspond to the image under some $A \in \text{GL}(2, \mathbb{Z})$ of the standard basis $\{ \mathfrak{t}_1, \mathfrak{t}_2\}$ of $\mathfrak{t}$.  (The proof by ``factoring'' the action is the same as in Proposition \ref{Smoothness}.)

Since the splitting is oriented, the vector $V_3$ extends uniquely as a unit trivialization of the line bundle along $C$.  We can use the complement subgroup $H$ to transport the vectors $\{V_1, V_2\}$ along $C$ to obtain an orthonormal frame $\{V_1, V_2, V_3\}$ for the bundle $N_xC$ such that $\{V_1,V_2\}$ span the 2-dimensional sub-bundle of the splitting of $NC$ and $V_3$ spans the line-bundle. (Note that the transportation of the vector $V_3$ along $C$ by $H$ agrees with the unique extension above.) Exponentiating this frame with respect to $\rho$ gives coordinates $\{\theta, x_1, x_2, x_3\}$ on a neighbourhood of $C$ such that along $S^1 \times \{0\}$, $\{ \frac{\partial}{\partial x_1}, \frac{\partial}{\partial x_2}\}$ span the 2 dimensional sub-bundle and $\frac{\partial}{\partial x_3}$ spans the line bundle.

Thus, in the coordinates on $N(C)$ induced by the frame and the coordinates on $T^2$ induced by a splitting $T^2 = H \times G_x$, $T^2$ acts as $(t_1, t_2) \cdot (\theta, r, \alpha, x_3) \mapsto (\theta + t_1, r, \alpha + t_2, x_3)$, where $r^2 = x_1^2 + x_2^2, \alpha = \arctan(x_2/x_1)$. The rest of the argument of Theorem \ref{hondaloc} goes through unmodified, using the equivariant version of the local Honda-Moser theorem as described above, applied to the given form $\omega$ and the form $\omega_A$ constructed in the above coordinates.
\end{proof}
\end{thm}

Given the reparametrization of $T^2$ via the splitting needed in the last theorem, the following fact describes how the moment map changes.

\begin{prop}
If $\Phi : \to \mathfrak{t}^*$ is a moment map for a Hamiltonian $T^2$ action $\sigma$, and $\mu: p \mapsto Ap + b \in \text{Aff}(2, \mathbb{Z})$ (i.e. $A \in GL(2,\mathbb{Z}), b \in \mathbb{R}^2$), then $\mu \circ \Phi$ is a moment map for the torus action $\sigma'( t, x) = \sigma(A^{-T}t, x)$.
\end{prop}
(Note that both actions have the same orbits.) 

Thus the moment map near a component $C$ of $Z_{\omega}$ has image $(A^{-T} \circ \Phi_0)(S^1 \times D^3) + b$, where $A \in GL(2, \mathbb{Z})$ is the integral transformation corresponding to the splitting used in the symplectomorphism to the standard fold, $\Phi_0$ is the moment map for the standard fold, and $b \in \mathbb{R}^2$.  Note that by definition, the image of the $p_1$ axis under $A^{-T}$ is independent of the choice of splitting (since it corresponds to $\mathfrak{g}_x^0$), so there is no ambiguity in the location of the fold. 

\subsection{Failure of convexity}\label{convexitysection}
Given the above, we have a completely canonical description of a neighborhood of any of the $T^2$ orbits.  It is natural to ask to what extent the theorems of the previous chapters apply.  In this section we consider convexity alone.  

First, local convexity still holds away from $Z_{\omega}$.  In some sense local convexity also holds near $Z_{\omega}$: the image of the moment map on a nice tubular neighborhood of a component $C$ is a convex subset of of $\mathbb{R}^2$.  In another sense, though, it fails, in that convex sets in the standard coordinates $(z, r^2)$ of the orbit space are not mapped to convex sets by the moment map $\Phi_0$.   

To see how global convexity fails, it's interesting to see how the Morse-Bott theory of the moments in the near-symplectic case differs from that in the symplectic case of Chapter 2.  

Let $\xi = (a,b) \in \mathfrak{t}$. For the standard fold, the moment $\Phi^{\xi}$ is given by
\begin{equation*}
\Phi^{\xi}( \alpha, x, y, z ) = a(z^2 - \frac{1}{2}(x^2+y^2)) + b( z(x^2+y^2))
\end{equation*}
and (suppressing the $\alpha$ direction by symmetry) its exterior derivative is 
\begin{equation*}
d\Phi^{\xi} = (2bz - a)x dx +  (2bz-a)y dy + (2az + b(x^2+y^2)) dz
\end{equation*}

For $(a,b) \neq (0,0)$, the critical set is described by two cases:
\begin{equation*}
\text{crit}(\Phi^{\xi}) = \begin{cases}
\{x=y=z=0\}& a \neq 0 \\
\{x=y=0\}& a = 0
\end{cases}
\end{equation*}
So in both cases, the critical set is a manifold. In these coordinates, the Hessian of $\Phi^{\xi}$ is
\begin{equation*}
H \Phi^{\xi} = 
\begin{pmatrix}
2bz-a &  0 & 2bx \\
 0 & 2bz-a & 2by \\
2bx & 2by & 2a 
\end{pmatrix} 
\end{equation*}
For the case $a \neq 0$, along the critical set $\{x=y=z=0\}$ (i.e. on C) we have 
\begin{equation} \label{bottcase}
H \Phi^{\xi} = 
\begin{pmatrix}
-a &  0 & 0 \\
 0 & -a & 0  \\
0 & 0 & 2a 
\end{pmatrix} 
\end{equation}
which is non-degenerate along the normal bundle (which is all of $D^3$), so $\Phi^{\xi}$ is  Morse-Bott.
For the case $a = 0$, along the critical set $\{x=y=0\}$ (i.e. on the line bundle) we have
\begin{equation}
H \Phi^{\xi} = 
\begin{pmatrix}
2bz &  0 & 0 \\
 0 & 2bz & 0  \\
0 & 0 & 0 
\end{pmatrix} 
\end{equation}
In this case the Hessian is non-degenerate in the normal directions (here being the x-y plane) everywhere except at $z=0$, and as $z$ crosses zero the stable manifold become unstable, so $\Phi^{\xi}$ is not Morse-Bott.  Note that the $z$-axis maps under $\Phi$ to the folded double $p_1$ axis, and that the change in sign corresponds to the inside normal to the image changing from down to up pointing as $z$ passes zero.

Recall, a neighborhood of each component of $Z_{\omega}$ is equivariantly symplectomorphic to the standard fold, up to an integral affine transformation of the moment map and torus.  So in general, the moments $\Phi^{\xi}$ fail to be Morse-Bott exactly for the finite set of $\xi$ for which $\xi \in \text{span}(A \frac{\partial}{\partial q_2})$, where $q_2$ is the standard coordinate on $T^2$ and and $A \in \text{GL}(2, \mathbb{Z})$ is the integral transformation of $T^2$ corresponding to the choice of splitting of $T^2$ to give the symplectomorphism to the standard fold near one of the components of $Z_{\omega}$.  Note however that this vector is, by construction of the symplectomorphism, always in the Lie algebra $\mathfrak{g}_x$ of the isotropy group $G_x$,  for $x \in C$ a component of $Z_{\omega}$, so is in fact invariantly defined. 

The failure of $\Phi^{\xi}$ to be Morse-Bott for this small set of $\xi$ might not be so bad if the other $\Phi^{\xi}$ behaved nicely. However, equation \ref{bottcase}, shows that, depending on the sign of $a$, some of the components $C$ of $Z_{\omega}$ will have stable manifolds of dimension 3, i.e. codimension 1.  The absence of such stable manifolds is exactly what was required to show uniqueness of local maxima. In effect, the stable manifolds of the components of $Z_{\omega}$ form hypersurfaces which separate the  stable manifolds of the local maxima.

\begin{figure}
\PSbox{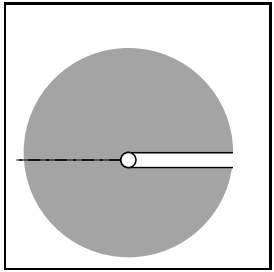}{1.3 in}{1.3 in}
\caption{The separating hypersurface for $b=0$, shown as the dotted line.}
\label{foldcut}
\end{figure}

We can draw this hypersurface easily in the case $b=0$, using the flat metric on $S^1 \times D^3$.  Then the gradient is (suppressing the $\alpha$ direction by symmetry)

\begin{equation}
\nabla \Phi^{\xi} = (-ax, -ay, 2az)
\end{equation}
so the stable manifold of $S^1 \times \{0\}$ is the set $\{z = 0\}$. This manifold is $T^2$ invariant, so can be described by its image under $\Phi_0$, which is the \emph{negative} $p_1$-axis.  This means that the separating hypersurface divides the manifold in a way that discounts local maxima/minima due to ``going around a fold''.  This is illustrated in Figure \ref{foldcut}.  Choosing $b \neq 0$ curves the line up or down but maintains its tangency to the $p_1$-axis at the origin.

For examples of near-symplectic toric 4-manifolds whose moment maps have non-convex images, see the construction in \ref{near-existence} and the final section.

\subsection{Classification theorems}
In this section I would like to extend the classification Delzant's theorem provides to near-symplectic manifolds.  The proof of uniqueness will follow with very little modification from that of Lerman-Tolman \cite{LT} presented in Chapter 2.  The existence proof will be a consequence of the existence of the local models developed so far, similar to the patching argument of Gay-Symington (\cite{Gay-S}), rather than the standard global existence proof for the Delzant theorem.

Because of possible overlaps in the image of the moment map $\Phi$, due both to folds and to the failure of global convexity, for the purpose of this section the analogue of the Delzant polytope $\Delta$ is the following abstract object:

\begin{defn}
A \emph{folded Delzant polygon}  is a triple $(B,F,\phi)$, where $B$ is a surface with corners, $F \subset (\partial B \setminus \{\text{corners}\})$ a discrete set, and $\phi: B \to \mathbb{R}^2$ a map that's a fold near $F$, an immersion on $B \setminus F$, takes edges to  line segments with rational slopes, and whose image satisfies the smoothness property of the Delzant theorem near the corners. 
\end{defn}

\begin{rem}
Here the ``smoothness property of the Delzant theorem'' is a strict version that requires the images of the corners to be the standard corner up to \emph{orientation-preserving} integral affine transformations; this implies local convexity at the corners.
\end{rem}

\begin{rem}
By  ``a fold near F'', we mean that there exist coordinates on $B$ near $F$ such that there $\phi$ is of the form $A \circ \phi_0 + b$, where $A \in \text{GL}(2,\mathbb{Z}), b \in \mathbb{R}^2$, and $\phi_0$ is the standard fold.  In the paragraph below, the surface $B$ will be the quotient $M/T$.
\end{rem}

For some examples of folded Delzant polygons, see Figures \ref{examp1-whole} and \ref{miscpolys}, which show the images $\phi(B) \subset \mathbb{R}^2$. Note that the immersions $\phi$ may fail to be 1-to-1 even away from the folds, and that the images may have an arbitrary number of ``holes''.

\subsubsection{Uniqueness}

By the canonical forms developed up to now, for any near-symplectic toric 4-manifold $(M, \omega)$ satisfying Condition \ref{invariantmetric} the orbit space $M / T^2$, the vanishing locus $Z_{\omega}$, and the map $\phi$ defined by $\Phi = \phi \circ \pi$ where $\pi: M \to M / T^2$ is the quotient, define a folded Delzant polygon $(M / T^2, \pi(Z_{\omega}), \phi)$.

We can state the uniqueness theorem, an analogue of Theorem \ref{delzant}, as follows.

\begin{thm}
Let $(M_1, \omega_1), (M_2, \omega_2)$ be two compact, connected, near-symplectic toric $4$-manifolds satisfying Condition \ref{invariantmetric} with moment maps $\Phi_1, \Phi_2$. Let $\pi_i: M_i \to B_i$ be the quotients to the orbit spaces, let $F_i = \pi_i (Z_{\omega_i})$  be the images of the vanishing locii, and define $\phi_i: B_i \to \mathbb{R}^2$ by $\Phi_i = \phi_i \circ \pi_i$. If there is a diffeomorphism $\psi: (B_1,F_1) \to (B_2,F_2)$ such that $\phi_2 \circ \psi = \phi_1$, then there exists a $T^2$-equivariant symplectomorphism $\Psi: M_1 \to M_2$ such that $\Phi_2 \circ \Psi = \Phi_1$.  
\begin{proof}
The proof is an adaptation of propositions \ref{deltadeterm} to \ref{ltsurj} as follows.  Propositions \ref{deltadeterm} and \ref{locsymp} are true at a point $x \in F$ by the model for the standard fold, and still true elsewhere.  

To check Proposition \ref{lambdatoh}, we need to show that the map $\Lambda$ is well-defined, i.e. that the vector field $X_f$ can be defined over $F$.  Write $\omega_A =  -2z(d \alpha dz + r dr d\theta) - r dr d\alpha + r^2 d\theta dz$.  Let $p_1 = z^2 -\frac{1}{2}r^2$, $p_2 = zr^2$ be coordinates on the base.  Then a general 1-form $\nu$ on the base can be written $\nu = a dp_1 + b dp_2$.  Solving $\iota(X)\omega_A = \pi^*(\nu)$, we obtain $X = -a \frac{\partial}{\partial \alpha} + b \frac{\partial}{\partial \theta}$. For $f \in C^{\infty}(B)$, $\nu = df$, have $a = \frac{\partial f}{\partial p_1}, b = \frac{\partial f}{\partial p_2}$. 

Here there is a technicality to worry about: if $f$ is just some smooth function on the base in the coordinates $(z,r^2)$, there's no guarantee that $a = \frac{\partial f}{\partial p_1}, b = \frac{\partial f}{\partial p_2}$ are well defined, or even bounded, at $F$.  However, if we require that $f$ is a smooth function on $\mathfrak{t}^*$ near $\phi(F)$, they are well-defined. So this fact requires a redefinition of the sheaf $\tilde{C^{\infty}}$ near $Z_{\omega}$.

Given the above, since $\frac{\partial}{\partial \theta} = -r \sin \theta \frac{\partial}{\partial x} + r \cos \theta \frac{\partial}{\partial y} $, $X_f$ extends over $F$ and there it is of the form $X_f = -a \frac{\partial}{\partial \alpha}$. Given this, the rest of the proposition generalizes by continuity.

The remaining propositions showing exactness of the sequence of sheaves go through unchanged.

The final step is to use the long exact sequence in sheaf cohomology to show that $H^1(B, \mathcal{H}) = 0$.  It is no longer necessarily the case that $B$ is contractible, or even simply-connected.  However, the relevant portion of the long exact sequence is:

\begin{equation*}
H^1(B, C^{\infty}) \to H^1(B, \mathcal{H}) \to H^2(B, \ell \times \mathbb{R})
\end{equation*}
The left hand term is still zero since $C^{\infty}$ is flabby.  I claim that $H^2(B, \ell \times \mathbb{R}) =0$: Since $B$ is a surface with non-empty boundary, it's homotopy-equivalent to a 1-complex, so it can be covered by contractible sets $\{U_{\alpha}\}$ such that no three intersect and each double intersection is contractible.  This means that the $\check{C}$ech cohomology $\check{H}^2( B, \ell \times \mathbb{R})$ of this cover is trivially zero, and the corresponding sheaf cohomology is as well (see, eg., R. O. Wells, \emph{Differential Analysis on Complex Manifolds}, p.64). 
\end{proof}
\end{thm}

\subsubsection{Existence}\label{near-existence}
The proof in this section follows ideas of Gay-Symington \cite{Gay-S} and Symington \cite{S1}.  

\begin{thm}
Let $(B,F,\phi)$ be a folded Delzant polygon.  Then there exists a near-symplectic toric 4-manifold $M$ with moment map $\Phi = \phi \circ \pi$, orbit space $B$, and vanishing locus $\pi^{-1}(F)$, where $\pi$ is the quotient by the $T^2$ action.
\begin{proof}
Consider the 4-manifold with boundary given by $\tilde{M} = B \times T^2$.  Set $\tilde{\Phi} = \phi \circ \tilde{\pi}$, where $\tilde{\pi}$ is projection on the first factor. Define a 2-form on $\tilde{M}$ by $\tilde{\omega} = \tilde{\Phi}^*(dp_1) \wedge dq_1 + \tilde{\Phi}^*(dp_2) \wedge dq_2$, where $(q_1,q_2)$ are standard coordinates on $T^2$ and $(p_1, p_2)$ are standard coordinates on $\mathbb{R}^2$.  Note that on $\pi^{-1}(\text{int}(B)) = \text{int}(\tilde{M})$, $\tilde{\omega}$ is symplectic, and the natural $T^2$ action given by multiplication in the second factor is Hamiltonian with moment map $\tilde{\Phi}$.

Now, define another manifold $M = \tilde{M} / \sim$ which is constructed by collapsing the fibres above the boundary $\partial B$ as follows:\\
1.  For $x \in \partial B \setminus ( F \cup \{\text{corners}\})$, collapse the $T^2$ fibre by taking the quotient by the $S^1$  subgroup generated by the 1-dimensional subspace of $\mathfrak{t}$ whose annihilator is parallel to the image of $d\phi _x$.\\
2. For $x \in F$, collapse the $T^2$ fibre by taking the quotient by the $S^1$  subgroup generated by the 1-dimensional subspace of $\mathfrak{t}$ whose annihilator is parallel to the image under $\phi$ of the edge that $x$ lies on.\\
3. For $x \in \{\text{corners}\}$, collapse the entire $T^2$ fibre.\\

The models developed in Chapter 2 guarantee that steps 1 and 3 can be done preserving the manifold structure, $T^2$ action, moment map, and symplectic form, such that the form on the quotient is symplectic there. The model developed at the beginning of this chapter guarantees that step 2 can be done preserving the manifold structure, $T^2$ action, moment map, and symplectic form, such that the symplectic form on the quotient vanishes over $F$. 
\end{proof}
\end{thm}

\subsection{Example(s)}

\begin{figure}
\PSbox{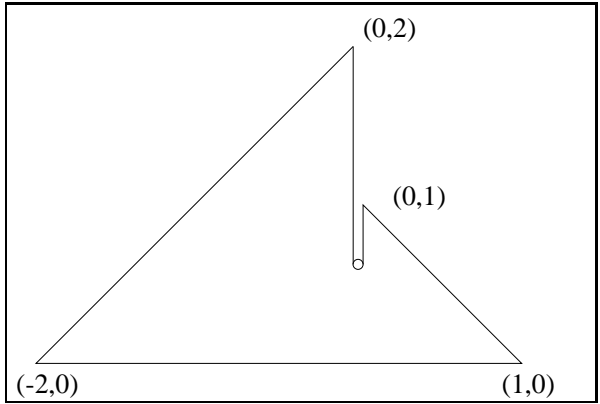}{3 in}{2 in}
\caption{The example, $\mathbb{C}P^2 \# \mathbb{C}P^2$}
\label{examp1-whole}
\end{figure}

In this section I analyze the near symplectic manifold corresponding to an example folded Delzant polygon in detail.  Consider the folded Delzant polygon illustrated in Figure \ref{examp1-whole}.  Here we show the image $\phi(B) \subset \mathfrak{t}^*$, and label the vertices with their coordinates.  $B$ is four-sided polygon with a single fold, indicated by the circle in the figure.  We draw two parallel lines to indicated the edges in the fold, though they really overlap. 

We can assemble the near-symplectic manifold given by Section \ref{near-existence} from a series of models, as follows.  First, consider the space $\mathbb{C}P^2 = (\mathbb{C}^3 \setminus \{0\}) / \sim$, where $\sim$ is multiplication by non-zero scalars, with the $T^2$ action induced from the $T^2$ action on $\mathbb{C}^3$ given by $(t_1,t_2) \cdot (z_1, z_2, z_3) = (e^{i t_1}z_1, e^{i t_2}z_2, z_3)$.  In homogeneous coordinates where $z_3 \neq 0$, define the moment map $\phi: \mathbb{C}P^2 \to \mathbb{R}^2$ to be 
\begin{equation}
\phi([z_1,z_2,1]) = ( \frac{|z_1|^2}{1 + |z_1|^2 + |z_2|^2}, \frac{|z_2|^2}{1+ |z_1|^2 + |z_2|^2})
\end{equation}

\begin{figure}
\PSbox{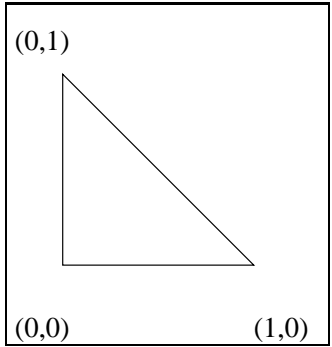}{1.6 in}{1.7 in}
\caption{The image $\phi(\mathbb{C}P^2) \subset \mathbb{R}^2$}
\label{cp2-std}
\end{figure}

Note that this moment map has image the right triangle shown in Figure \ref{cp2-std}, with vertices at $\{(0,1), (1,0), (0,0)\}$, and in the angle coordinates induced by the torus action, the collapsing of the $T^2$ fibres over the edges is as described in Section \ref{near-existence}.

We calculate, for $j,k \in \{1,2\}, j \neq k$, $(p_1,p_2)$ coordinates on $\mathbb{R}^2$, 

\begin{equation}
\phi^*(dp_j) = \frac{1}{(|z_1|^2 + |z_2|^2 + 1)^2} ( (|z_k|^2 +1)(z_j d \bar{z_j} + \bar{z_j} dz_j) - |z_j|^2(z_k d \bar{z_k} + \bar{z_k}dz_k) )
\end{equation} 

Setting $z_j = r_j e^{2 \pi i \theta_j}$, we have 

\begin{equation}
d \theta_j = -\frac{i}{2} \frac{\bar{z_j} dz_j - z_j d \bar{z_j}}{|z_j|^2}
\end{equation}

So the symplectic form corresponding to this moment map is 

\begin{equation}
\begin{split}
d\theta_1 \wedge \phi^*(dp_1) + d \theta_2 \wedge \phi^*(dp_2) =&  \frac{-i}{(|z_1|^2 + |z_2|^2 + 1)^2}[  (|z_2|+1)dz_1 \wedge d\bar{z_1} + (|z_1|^2+1)dz_2 \wedge d\bar{z_2}] \\
& +  \frac{-i}{(|z_1|^2 + |z_2|^2 + 1)^2} [z_1 \bar{z_2} d \bar{z_1} \wedge dz_2 - \bar{z_1} z_2 dz_1 \wedge d \bar{z_2}]               
\end{split}
\end{equation}

This is the Fubini-Study form $\omega_{FS}$ on $\mathbb{C}P^2$, so we don't have to check that it extends symplectically to $\cup_{j=1}^3 \{z_j=0\}$. 

\begin{figure}
\PSbox{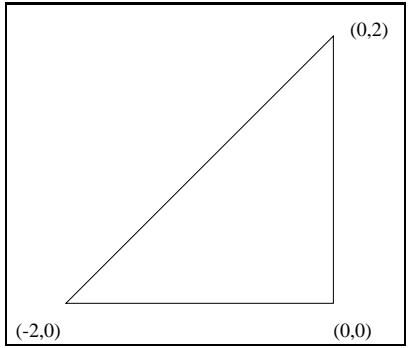}{2 in}{1.7 in}
\caption{The image $\tilde{\phi}(\mathbb{C}P^2) \subset \mathbb{R}^2$}
\label{cp2-2}
\end{figure}

Similarly, consider the space $\mathbb{C}P^2$ with the $T^2$ action induced by the action on $\mathbb{C}^3$ given by $(t_1,t_2) \cdot (z_1, z_2, z_3) = (e^{-it_1}z_1, e^{it2}z_2, z_3)$, and moment map 
\begin{equation}
\tilde{\phi}([z_1,z_2,1]) = ( \frac{-2|z_1|^2}{1 + |z_1|^2 + |z_2|^2}, \frac{2|z_2|^2}{1+ |z_1|^2 + |z_2|^2})
\end{equation}

Note that this moment map has image the right triangle shown in Figure \ref{cp2-2}, with vertices at $\{(0,2), (-2,0), (0,0)\}$.

The new $T^2$ action gives a new angular coordinate $\tilde{\theta_1} = -\theta_1$, where $z_1 = r_1 e^{2 \pi i \theta_1}$. Combined with the new moment map $\tilde{\phi}$, the sign changes cancel and this gives the induced symplectic form

\begin{equation}
d\tilde{\theta_1} \wedge \tilde{\phi}^*(dp_1) + d \theta_2 \wedge \tilde{\phi}^*(dp_2) = 2 \omega_{FS}
\end{equation}

which is twice the Fubini-Study form, so again it extends symplectically to $\cup_{j=1}^3 \{z_j = 0\}$.

The last model we need is the space $\mathbb{C} \times \mathbb{R} \times S^1$.  Let $(r, \theta)$ be polar coordinates on $\mathbb{C}$, $x$ a coordinate on $\mathbb{R}$, and $\alpha$ a coordinate on $S^1$.   Let $T^2$ act by $(t_1, t_2)\cdot(r, \theta, x, \alpha) = (r, \theta+t_2, x, \alpha+t_1)$. The moment map $\phi(r, \theta, x, \alpha) = (x, r^2)$ induces the symplectic form $d \theta_1 \wedge \phi^*(dp_1) +  d \theta_2 \wedge \phi^*(dp_2) = d \alpha \wedge  dx + d \theta \wedge 2r dr$ on $\mathbb{C} \times \mathbb{R} \times S^1$.

We construct the example space as follows.  We remove from $(\mathbb{C}P^2, \omega_{FS})$ the ball corresponding to $\phi_1 + \phi_2 < 5/8$, i.e. set 
\begin{equation}
M_A = (\mathbb{C}P^2, \omega_{FS}) \setminus \{  \frac{|z_1|^2}{1 + |z_1|^2 + |z_2|^2}+ \frac{|z_2|^2}{1+ |z_1|^2 + |z_2|^2} < 5/8 \}
\end{equation}

Similarly, remove a ball from $(\mathbb{C}P^2, 2\omega_{FS})$ by setting

\begin{equation}
M_B = (\mathbb{C}P^2, 2\omega_{FS}) \setminus \{  \frac{2|z_1|^2}{1 + |z_1|^2 + |z_2|^2}+ \frac{2|z_2|^2}{1+ |z_1|^2 + |z_2|^2} < 5/8 \}
\end{equation}

We restrict the model $\mathbb{C} \times \mathbb{R} \times S^1$ to the set $r^2 < 3/8$, $r^2 + x < 3/4$, $r^2 - x < 3/4$ to obtain the manifold $M_C$.  

Finally, we rotate the standard fold by ninety degrees, shift up by $1/2$, and restrict to the set $p_2 > 2/8,  p_1 + p_2 < 3/4$, $-p_1 + p_2 < 3/4$ to obtain the manifold $M_D$. 

\begin{figure}
\PSbox{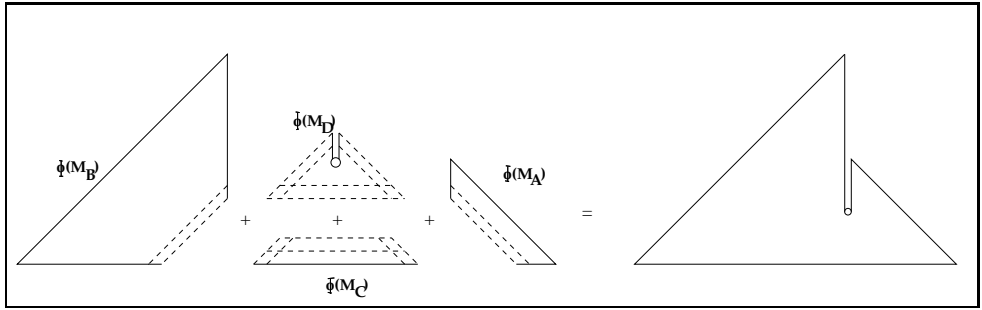}{5in}{1.5 in}
\caption{Patching the models}
\label{examp1-parts}
\end{figure}

The last step is identifying the four manifolds over the strips indicated in Figure \ref{examp1-parts}. Solid lines indicate edges over which the  $T^2$-fibres are partially or completely (at corners) collapsed, while dashed lines indicate the open sets where patching will take place.  We patch $M_A$, $M_C$, and $M_D$ on  $5/8 < p_1 + p_2 < 3/4$, $M_B$, $M_C$, and $M_D$ on $5/8 < -p_1 + p_2 < 3/4$, and $M_C$ and $M_D$ on $1/4< p_2 < 3/8$, via the coordinates $(p_1, \theta_1, p_2, \theta_2)$ in which all four manifolds have symplectic form $d \theta_1 \wedge dp_1 +  d \theta_2 \wedge dp_2$. Since all the models are symplectic, the fibres are collapsed the same way over the edges, and the identification is a symplectomorphism on the open dense set where the fibres are $T^2$, the identification is a symplectomorphism everywhere.

\begin{figure}
\PSbox{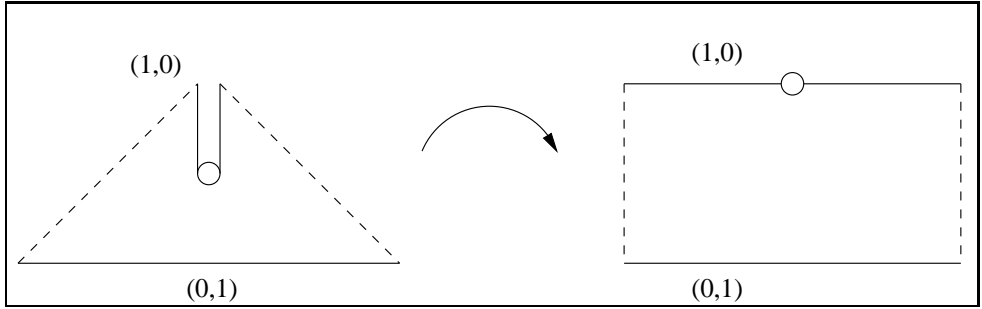}{5in}{1.5 in}
\caption{The patch, $S^3 \times [0,1]$, with collapsing directions labeled}
\label{patch}
\end{figure}

We claim that the resulting manifold is $\mathbb{C}P^2 \# \mathbb{C}P^2$:   We've connected two copies of $\mathbb{C}P^2 \setminus \{\text{ball}\}$  using a patch whose image under the moment map is given on the left side of Figure \ref{patch}. (Here each edge is labelled with the generator of the subgroup which is collapsed in the fibre.) The patch is diffeomorphic to a trival $T^2$-fibration over $[0,1] \times [0,1]$, on which we collapse the fibres over $[0,1] \times \{0\}$ along the first $S^1$ factor and we collapse the fibres over $[0,1] \times \{1\}$ along the second $S^1$ factor, as shown on the right of Figure \ref{patch}. This fibration, finally, is diffeomorphic to $S^3 \times [0,1]$, as can be seen by considering the model $S^3 \times [0,1] = \{ (z_1, z_2, x) \in \mathbb{C}^2 \times \mathbb{R} | |z_1|^2 + |z_2|^2 = 1, x \in [0,1]\}$ with the $T^2$-action $(t_1, t_2)\cdot(z_1, z_2,x) = (e^{2 \pi i t_1}z_1, e^{2 \pi i t_2}z_2,x)$ and the projection $(z_1, z_2,x) \mapsto (x, |z_1|^2)$.  The patching corresponds to that for the connected sum.

\begin{figure}
\PSbox{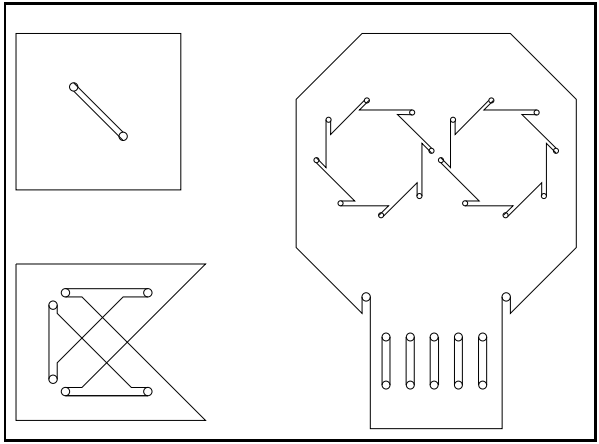}{3 in}{2.2 in}
\caption{Assorted folded Delzant polygons. The left two come from \cite{Gay-S}.}
\label{miscpolys}
\end{figure}

The construction in Section \ref{near-existence} guarantees the existence of near-symplectic toric 4-manifolds of a wide variety which may be less familiar.  Some folded Delzant polygons giving rise to such manifolds are given in Figure \ref{miscpolys}.  

\appendix
\section{Construction of invariant metrics}\label{invariantmetricproof}
Given a near-symplectic G-invariant form $\omega$, the following argument, explained to me by Denis Auroux and based on the non-equivariant argument in \cite{Au}, p.63, guarantees the existence of a G-invariant Riemannian metric $\rho_{\omega}$ with respect to which $\omega$ is self-dual.  Note that transverse vanishing is independent of the metric, so this guarantees that Condition \ref{invariantmetric} is always satisfied.

Consider a vector space $V \cong \mathbb{R}^4$ with a fixed postive-definite inner product $\rho_0$ and a given orientation. The Hodge-* operator corresponding to $\rho_0$ and the Hodge inner-product on $\wedge^2 V$, $\langle \eta, \nu \rangle = \eta \wedge * \nu$, induce an orthogonal splitting of $\wedge^2 V = \wedge^2_{+,0} \oplus \wedge^2_{-,0}$ into the $\rho_0$-self-dual and anti-self-dual 2-forms.  For any other positive-definite inner product $\rho_i$ on $V$, its corresponding space of self-dual 2-forms $\wedge^2_{+,i}$ is a 3-plane in $\wedge^2 V$ on which the wedge-product restricts to a positive-definite bilinear form.  Any such 3-plane can be written uniquely as the graph $P = \{ \alpha + L_i(\alpha), \alpha \in \Lambda^2_{+,0}\}$ of a linear map $L_i: \wedge^2_{+,0} \to \wedge^2_{-,0}$ with operator norm less than 1.  Conversely, any such linear map $L_i$ defines a positive definite inner product $\rho_i$ on V up to scaling, by specifying its space of self-dual 2-forms.  Note that this space of maps is convex.

On the $4$-manifold $M$ with a $G$-action and a $G$-invariant near-symplectic form $\omega$, choose any $G$-invariant Riemannian metric $\rho_0$.  This induces a $G$-invariant splitting of the bundle of $2$-forms $\wedge^2 (T^*M) = \wedge^2_{+,0}(T^*M) \oplus \wedge^2_{-,0}(T^*M)$. Because $\omega$ is near-symplectic, it's  self-dual with respect to some other (non-invariant) Riemannian metric $\rho_1$.  By the discussion above, the bundle of $\rho_1$'s self-dual forms is the graph of a section $L_1 \in \text{Hom}(\wedge^2_{+,0}(T^*M), \wedge^2_{-,0}(T^*M))$ with pointwise operator norm less than 1, and $\omega$ is in this graph.  Since the splitting $\wedge^2 (T^*M) = \wedge^2_{+,0}(T^*M) \oplus \wedge^2_{-,0}(T^*M)$ is $G$-invariant, we can average $L_1$ over $G$, using the convexity above, to obtain a $G$-equivariant section $\tilde{L_1} \in \text{Hom}(\wedge^2_{+,0}(T^*M), \wedge^2_{-,0}(T^*M))$, which still has pointwise operator norm less than 1, and has $G$-invariant graph.  Since $\omega$ is $G$-invariant, $\omega$ is still in this graph. The graph of $\tilde{L_1}$ defines a conformal class of metrics having it as their bundle of self-dual 2-forms.  Take one such metric $\tilde{\rho_1}$.  Since $\tilde{L_1}$ has $G$-invariant graph, averaging $\tilde{\rho_1}$ over $G$  preserves the conformal class and produces a $G$-invariant Riemannian metric $\rho_{\omega}$ on $M$ having $\omega$ as a self-dual $2$-form.

\section{Hamiltonian $S^1$ actions}

One obvious direction for further work is to see whether a generalization of Karshon's classification of 4-manifolds with Hamiltonian $S^1$-actions is true in the near-symplectic case, and in particular whether her result that ``isolated fixed points implies toric variety'' is true.  

Honda's local model $\omega_B$ for the unoriented splitting provides a local counterexample to the toric claim, in that it has a Hamiltonian $S^1$ action but no $T^2$ action.  I do not know if there is a compact example containing an unoriented splitting, or if the presence of unoriented splittings is the only obstruction.
In any case, the analysis using the equivariant Honda-Moser theorems used in Chapter 4 provides a first step to generalizing Karshon's results, by describing local models for an $S^1$ action in a neighbourhood of $Z_{\omega}$.  In particular, I claim the following:

 Let $x \in C$, a component of $Z_{\omega}$.

(1) If $x$ has trivial stabilizer, then a neighbourhood of $C$ is equivariantly symplectomorphic to the standard fold, with Hamiltonian $S^1$ action given by the moment $\Phi_1 = z^2 - \frac{1}{2}r^2$.

(2) If $x$ has $S^1$ stabilizer, then a neighbourhood of $C$ is equivariantly symplectomorphic to the standard fold, with Hamiltonian $S^1$ action given by the moment $\Phi_2 = zr^2$. 

(3) If $x$ has stabilizer $\mathbb{Z}_k \subset S^1$, and the splitting is oriented, a neighbourhood of $C$ is equivariantly symplectomorphic to the standard fold with Hamiltonian $S^1$ action given by the moment $k \Phi_1 + \Phi_2$.

(4) If the splitting is unoriented, $x$ must have stabilizer $\mathbb{Z}_2$, and a neighborhood of $C$ is equivariantly symplectomorphic to the unoriented model $\omega_B$ with $S^1$ action given by rotation in $\alpha$ (i.e. along $C$) and moment map $z^2 - \frac{1}{2}r^2$, which is well-defined on the quotient.

The proofs are similar to the proof of Theorem \ref{canfoldthm}. Given these local models, their local Morse theory is as described in Section \ref{convexitysection}.  This is very different from the symplectic case as analyzed by Karshon and will require a different analysis.  In particular, the existence of multiple local maxima/minima and the separating hypersurfaces require special attention.


\end{document}